\documentclass[11pt,letterpaper]{article}
\usepackage[utf8]{inputenc}
\usepackage{amsmath}
\usepackage{amsfonts}
\usepackage{amssymb}
\usepackage{amsthm}
\usepackage{graphicx}
\usepackage{xcolor}
\usepackage{url}
\usepackage[left=2cm,right=2cm,top=2cm,bottom=2cm]{geometry}

\theoremstyle{definition}
\newtheorem{theorem}{Theorem}

\newtheorem{lemma}{Lemma}

\newcommand{\overlines}[1]{#1}

\newcommand{\POP}{\textbf{\textit{POP}}}

\author{Will Traves and David Wehlau}
\title{Ten Points on a Cubic}

\begin{document}
\maketitle
\begin{abstract}
    In 1639 the 16-year old Blaise Pascal found a way to determine if 6 points lie on a conic using a straightedge. We develop a method that uses a straightedge to check whether 10 points lie on a plane cubic curve. 
\end{abstract}

\setlength{\parskip}{1em}

Books I through VI of Euclid's \emph{Elements} treat ruler and compass constructions, which remain a mainstay of high school geometry in America today. Two classical problems that cannot be solved by ruler and compass are to double the cube -- to construct a line segment of length $\sqrt[3]{2}$ given a segment of length 1 -- and to trisect a general angle. Today these impossibility results are usually proven using Galois theory, but the first proof dates to 1837, just five years after Evariste Galois was killed in a duel. That proof \cite{Wantzel}, given by the 23-year old mathematician Laurent Wantzel, was ignored and forgotten for over 80 years. There are many variants of the constructibility problems. In some we use a rusty compass which only has one setting. In another we replace the compass by the ring left from a coffee cup: we have no compass but are given a single circle. Our favorite is to toss the compass away completely and just consider straightedge-only constructions! Using a straightedge, we are only allowed to draw lines between known points and construct points by intersecting lines. This requires us to produce incidence relations -- three collinear points or three concurrent lines -- to characterize geometric properties. It may be surprising that even with this reduced material there are many beautiful results. 

The 16-year old Blaise Pascal found a nice incidence result characterizing points on a conic in 1639. If six points lie on a conic then connecting the points with line segments forming a path gives a (degenerate) hexagon whose three pairs of opposite sides extend to meet in three \emph{collinear} points. The Pascal line through these three auxiliary points is depicted in Figure 1 as a dotted line. There are $6!=720$ ways to reorder the points but since rotations and reflections of the hexagon do not produce a new picture, we can draw just 60 hexagons, each giving rise to a Pascal line. The reader is invited to choose a different hexagon and draw the Pascal line. This arrangement of 60 lines is known as Pascal's Hexagrammum Mysticum and has been studied by many important geometers, including the Reverend T.P. Kirkman, Arthur Cayley, Jakob Steiner, Julius Pl\"ucker and George Salmon. See Conway and Ryba's papers \cite{CRy1, CRy2} for details, including a reference to the hand-drawn diagram of all 60 lines due to Anne and Elizabeth Linton \cite{Linton}, twin sisters who completed doctoral studies together at the University of Pennsylvania in 1921.

\begin{figure}[h!t]
\begin{center}
\includegraphics[scale=0.2]{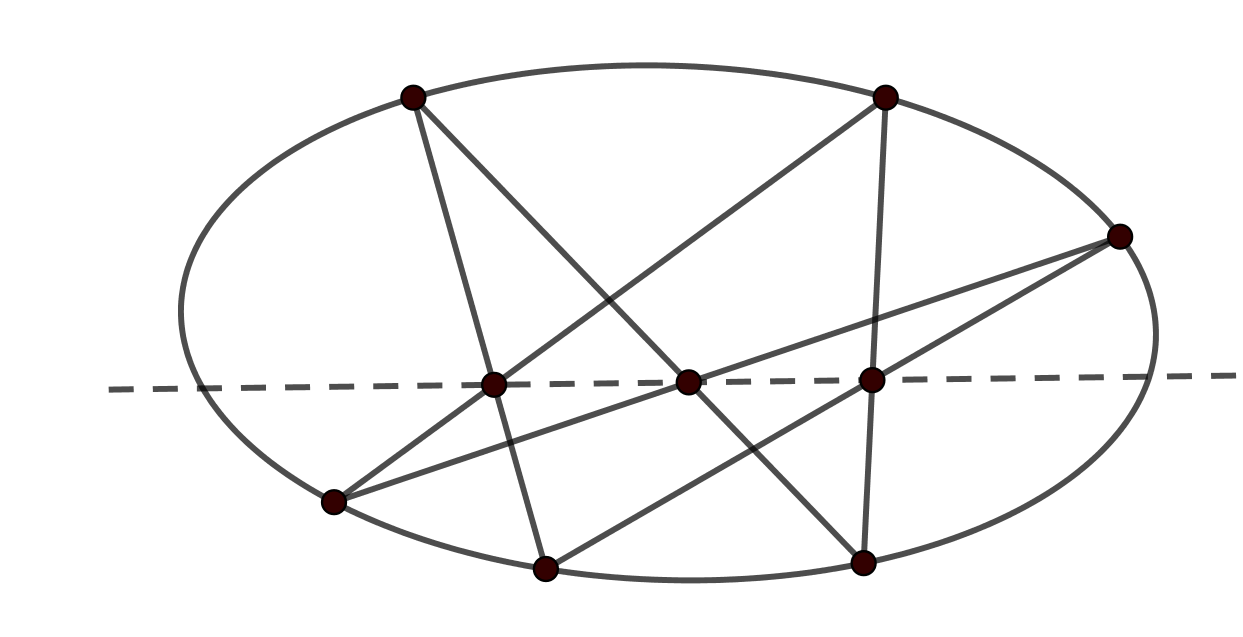}
\end{center}
\vspace{-0.3in}
\caption{An incidence result relating six points on a conic to 3 points on a line.}
\label{PBM}
\end{figure}

 The English mathematicians William Braikenridge and Colin Maclaurin established the converse to Pascal's Theorem almost a hundred years after Pascal's discovery: if three lines meet another set of three lines in 9 points with three of the points lying on yet another line then the remaining six points lie on a conic. Today we know these theorems are part of a family of such results described by the Cayley-Bacharach Theorem.  

\begin{figure}[h!t]
\begin{center}
\includegraphics[scale=0.2]{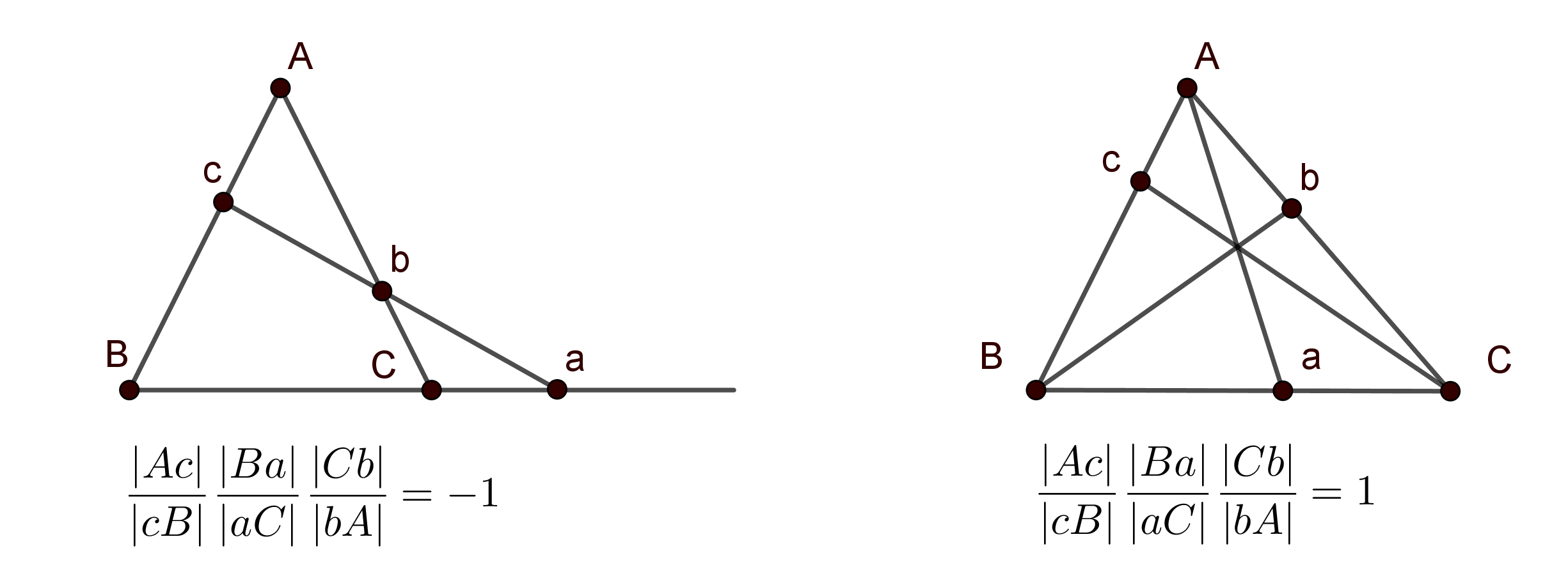}
\end{center}
\vspace{-0.3in}
\caption{Menelaus's Theorem (left) and Ceva's Theorem (right).}
\label{MC}
\end{figure}

While a straightedge construction cannot use angles or distances directly, many incidence conditions are equivalent to angle or distance constraints. The best known of these are Menelaus's and Ceva's Theorems.  The first theorem, due to Menelaus of Alexandria (70-140 A.D.) says that three points on the three (extended) edges of a triangle are collinear precisely when the product of three oriented length ratios is -1, as illustrated in Figure \ref{MC} (left). Ceva's Theorem was first proven by Yusuf Al-Mu'taman ibn H\"ud, an eleventh-century king of Saragossa in present-day Spain, and later proven and popularized by Giovanni Ceva in 1678. Ceva's Theorem dualizes Menelaus's Theorem, interchanging lines and points: the sides of a triangle are cut by three concurrent lines that pass through the corresponding opposite vertex precisely when the product of the three oriented length ratios is 1, as illustrated in Figure \ref{MC} (right).   

\begin{figure}[h!t]
\begin{center}
\includegraphics[scale=0.2]{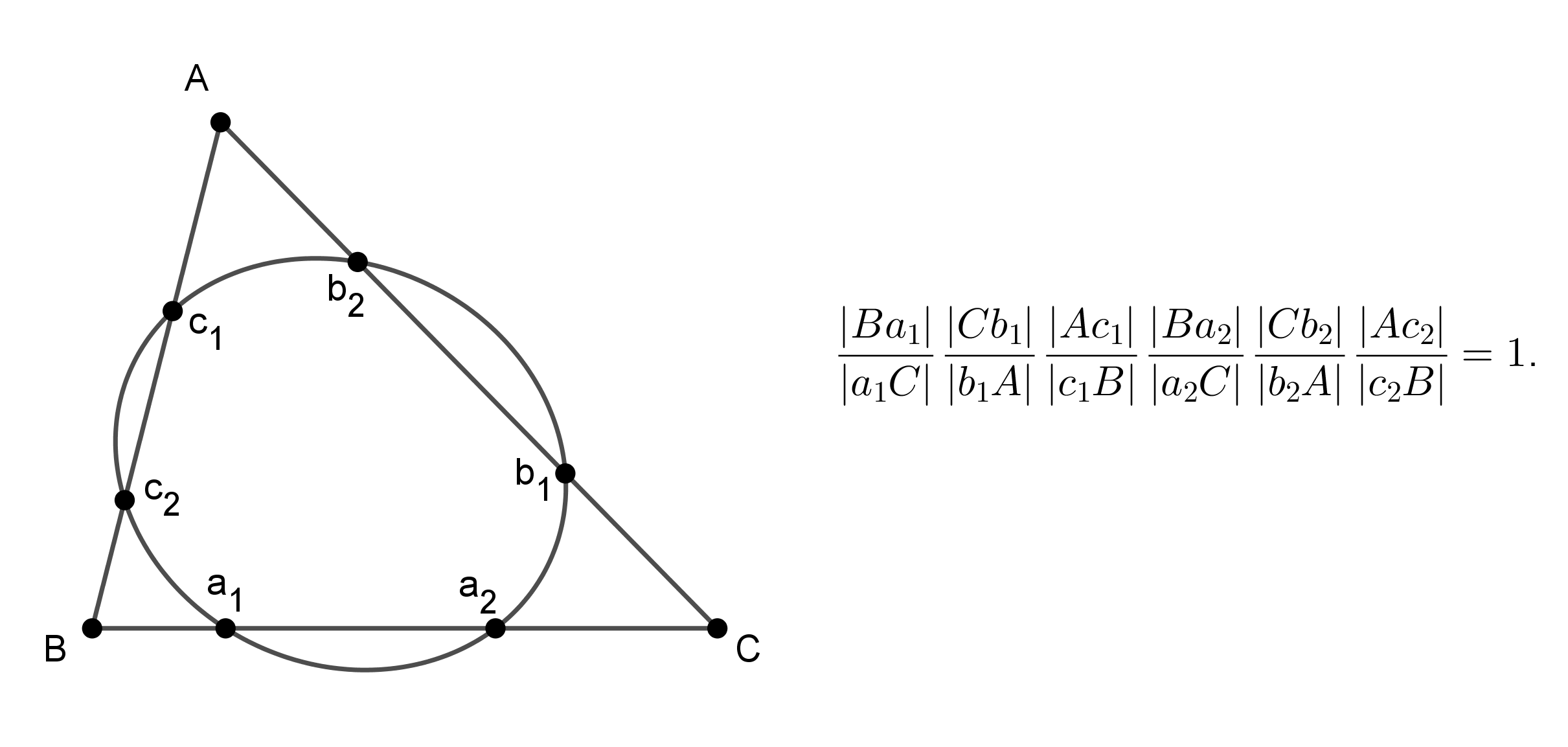}
\end{center}
\vspace{-0.3in}
\caption{Carnot's Theorem: labeling the six points (left) and the product condition (right).}
\label{CT}
\end{figure}

Lazare Carnot is best known as the ``Organizer of Victory" for the French Revolutionary Army at the end of the 18$^\text{th}$ century but after being exiled for his revolutionary activities he retired to write about mathematics and military tactics. Like Pascal's Theorem, Carnot's Theorem characterizes when six points lie on a conic, but the result involves products of distance ratios like Menelaus's and Ceva's theorems. Carnot drew three lines through pairs of the six points, producing the triangle $ABC$. Labeling the points $a_1$, $a_2$, $b_1$, $b_2$, $c_1$ and $c_2$ as in Figure~\ref{CT}, Carnot observed that the six points lie on a conic precisely when a product of six distance ratios equals 1. 

A natural next step is to look for incidence theorems involving cubic curves. One consequence of the Cayley-Bacharach Theorem is the \emph{Eight Implies Nine Theorem}: Given any eight points, there is a ninth special point so that every cubic through the eight given points must pass through the ninth point. Using a straightedge to construct the ninth point from the eight given points was a problem considered by several people in the 19th century, including A. S. Hart \cite{Hart}, Thomas Wheedle \cite{Wheedle}, Michel Chasles \cite{Chasles}, and Arthur Cayley \cite{Cayley}.  Recently, Qingchun Ren, J\"{u}rgen Richter-Gebert, and Bernd Sturmfels \cite{RRGS} took a modern approach to the problem. 
Though the Cayley-Bacharach Theorem plays a vital role in our work, we do not even state it here. The curious reader is referred to the excellent survey paper by Eisenbud, Green and Harris \cite{EGH}, who trace the history of the theorem and prove several versions of the theorem.  

Our main result gives a straightedge construction that checks whether 10 general points lie on a cubic. We give an incidence result reminiscent of the results of Pascal, Braikenridge, and Maclaurin: ten points lie on a cubic precisely when three constructed points are collinear. The key step is to realize the 10 points as a subset of 16 points that are the intersection of two degree-4 curves, each the union of two conics. The Cayley-Bacharach Theorem implies that the 10 points lie on a cubic precisely when the 6 remaining points lie on a conic, which we test using Carnot's Theorem. A serious complication is that when using a straightedge we can only find two of the 6 extra points together with two lines that contain the remaining four points. Remarkably, we are able to overcome this difficulty using important methods from geometry, such as the power of a point and cross ratios. These tools allow us to use Carnot's Theorem to check whether the six residual points lie on a conic, reducing the computation to an instance of Menelaus's Theorem. 

The next section describes the geometric context and introduces some mathematical tools. Section~\ref{GCsection} deals with constructive problems in synthetic geometry involving conics and lines, results that will be needed in our construction. Section~\ref{Section:10} gives the construction to check whether 10 points lie on a cubic. We give a second, computer-aided proof, that the construction works in Section~\ref{sec:bp}. The last section contains some pointers to the literature and several exercises for the interested reader.

\begin{section}{Context}

\noindent {\bf The plane.} We will be working with lines and points in the plane, but even this simple statement needs some clarification. It will be convenient to allow complex coordinates in many of our proofs so we'll work with $\mathbb{C}^2$ rather than $\mathbb{R}^2$. This means that our straightedge is a \emph{complex} straightedge: it allows us to draw the complex line $\{(x_0 + x_1t, \; y_0+y_1t): t \in \mathbb{C}\}$ joining two complex points $(x_0,y_0)$ and $(x_1,y_1) $ in $\mathbb{C}^2$ and to intersect two complex lines. Fortunately, if all of the geometric objects used as inputs to our constructions are defined over a subfield (like $\mathbb{R}$, $\mathbb{Q}$, or $\mathbb{Q}(\sqrt{2})$) of $\mathbb{C}$ then our output will also be defined over the same field. 

Moreover we work in the \emph{projective} plane $\mathbb{P}^2$, a copy of the usual plane together with \emph{points at infinity}. More generally, the points in the projective space $\mathbb{P}^n$ are modeled as one-dimensional subspaces of $\mathbb{C}^{n+1}$, lines through the origin $\mathbf{0} \in \mathbb{C}^{n+1}$. Each point in $\mathbb{P}^n$ is defined by its \emph{homogeneous} coordinates: the subspace with basis given by the vector $\langle x_1, x_2, \ldots,x_{n+1} \rangle \neq \mathbf{0}$ is denoted $[x_1:x_2:\ldots:x_{n+1}]$ (or just $[x:y:z]$ for points in the projective plane), where the square brackets remind us that a point is really an equivalence class of all possible basis vectors for the subspace. Since the basis vectors for a fixed one-dimensional subspace are nonzero scalar multiples of one another, this forces the equality $[\lambda x_0: \lambda x_1: \ldots, \lambda x_n] = [x_1: x_2: \ldots : x_{n+1}]$ for all non-zero scalars $\lambda \in \mathbb{C}$. The colon in the notation indicates that it is the ratios of these coordinates that determine the point in the projective space; the actual values of the coordinates are not so important since they can always be multiplied by a common scalar. It is customary to identify points $[x_0:x_1:\ldots: x_{n+1}]$ satisfying $x_{n+1} \neq 0$ with points in the usual $\mathbb{C}^n$: the point $[x_1:x_2:\ldots:x_{n+1}] = [x_1/x_{n+1}: x_2/x_{n+1}: \cdots: 1]$ is identified with $(x_0/x_{n+1},x_1/x_{n+1},\ldots,x_{n}/x_{n+1})$. The remaining points, all of the form $[x_1:x_2:\ldots:x_n:0]$, are thought of as points \emph{at infinity}. To get a sense of how the regular points connect with the points at infinity, we can traverse a line starting at
$\mathbf{b} = (b_1,\ldots,b_n) \in \mathbb{C}^n$, moving in the direction $\mathbf{m} = \langle m_1, \ldots, m_n \rangle \in \mathbb{C}^n$, and take a limit: 

$$\begin{array}{lll}  \displaystyle\lim_{t\rightarrow \pm \infty} (m_1t+b_1,\ldots,m_nt+b_n) & = & \displaystyle\lim_{t\rightarrow \pm \infty} [m_1t+b_1:\ldots: m_nt+b_n: 1] \\ &  = & \displaystyle\lim_{t\rightarrow \pm \infty} \displaystyle\left[m_1+\frac{b_1}{t}:\ldots:m_n+\frac{b_n}{t}:\frac{1}{t}\right] \\ &  = & [m_1:\ldots:m_n:0]. \end{array}$$
The limiting point at infinity $[m_1:\ldots:m_n:0]$ corresponds to the direction vector of the line, irrespective of the starting point, the speed, or even which way we move along the line.  

Just as points in $\mathbb{P}^2$ are modeled by one-dimensional subspaces of $\mathbb{C}^3$, lines in the projective plane are modeled by two-dimensional subspaces of $\mathbb{C}^3$. Such a subspace is a plane through the origin and is completely determined by its normal vector $\langle a,b,c \rangle \neq \langle 0,0,0 \rangle$: $L = \{[x:y:z] \in \mathbb{P}^2: \; ax+by+cz = 0\}$. Note that all the points at infinity lie on the line $z=0$, the \emph{line at infinity}. The line through distinct points $P$ and $Q$ has normal vector given by the cross product $P \times Q$. The cross product $v \times w$ also gives the homogeneous coordinates of the point of intersection of two distinct lines with normal vectors $v$ and $w$. In particular, two different lines in $\mathbb{P}^2$ always meet in a point: the point lies at infinity if the lines are parallel. We will use an \emph{infinitely long} straightedge. That is, we can construct the common point at infinity of two parallel lines. If the line at infinity can be constructed from our given input data then we can also construct the intersection point of any given line with the line at infinity.  The collection of lines in $\mathbb{P}^2$ is itself a 2 dimensional projective space, called the dual projective plane. The line with equation $ax+by+cz=0$ is identified with the point $[a:b:c]$. Note that this identification is well-defined since any other equation of the line is a nonzero scalar multiple $\lambda ax+\lambda by+\lambda cz=0$ of the original equation and its associated point $[\lambda a: \lambda b: \lambda c]$ equals $[a:b:c]$.

\noindent{\bf Circles as special conics.} Given a curve defined by a polynomial equation $g(x,y)=0$ in the usual plane, there is a standard way to extend the curve to all of $\mathbb{P}^2$. If the polynomial $g$ has degree $d$, we \emph{homogenize} $g$ by multiplying each term in the polynomial by a power of $z$ to ensure that it too has total degree $d$. The homogenization $G(x,y,z)$ vanishes at all points in the regular plane where $g$ vanished since $G(x,y,1)=g(x,y)$. While we cannot talk about the value of $G$ at a point $[x:y:z]$ because the value of $G$ changes when we scale the point, $G(\lambda x, \lambda y, \lambda z) = \lambda^d G(x,y,z)$, the \emph{vanishing} of $G(x,y,z)$ is well-defined. For instance, homogenizing the line cut out by $y-(mx+b)=0$ gives the projective line $y - (mx+bz) = 0$, which meets the line at infinity $z=0$ at $[1:m:0]$. As a second example, a general conic is defined by the vanishing of a degree-2 polynomial $$g(x,y) = ax^2 + bxy + cy^2 + dx + ey + f$$ and its extension to the full projective plane is cut out by setting the homogenization $$G(x,y,z) = ax^2 + bxy + cy^2 + dxz + eyz + fz^2$$ equal to 0. Note that imposing the requirement that a given point $P[x_0:y_0:z_0]$ lies on the conic $G(x,y,z)=0$ imposes a linear condition on the coefficients $a$, $b$, $\ldots$, $f$ of the polynomial $G$. Imposing five such conditions gives a $5 \times 6$ system of linear equations, which always has a nonzero solution by the Rank-Nullity Theorem, so there is a conic through any five points. However, the $6 \times 6$ system of linear equations imposed by six general points has a nonzero solution precisely when the six points lie on a conic. The extension of these ideas to degree-$d$ polynomials show there is a degree-$d$ curve through any $\binom{d+2}{2}-1$ points but $\binom{d+2}{2}$ points need to be in special position to lie on such a curve. In particular, 10 points need to be in special position in order to lie on a cubic curve. Our goal in this paper is to give a criterion describing these special positions in terms of incidence relations. 

One feature of working with curves in the projective plane is that intersections are easy to describe. Isaac Newton first observed in the Principia that two projective plane curves of degrees $d$ and $e$ meeting in finitely many points actually meet in precisely $de$ points, counted appropriately. Over a hundred years later, \'{E}tienne B\'ezout \cite{Bezout} generalized this observation to geometric objects in higher dimensional projective spaces and the result is generally known as B\'{e}zout's Theorem. For instance, as we've already seen, any two lines (degree-1 curves) meet in one point, though that point might lie at infinity. Lines meet conics in two points, but we need to count the intersection properly: a tangent point is counted twice (as in Figure~\ref{f2} (middle)) and the intersection points may be complex (and so not visible in the real drawing, as in Figure~\ref{f2} (right)). See Fulton \cite{Fulton} for details. 

\begin{figure}[h!t]
\begin{center}
\includegraphics[scale=0.4]{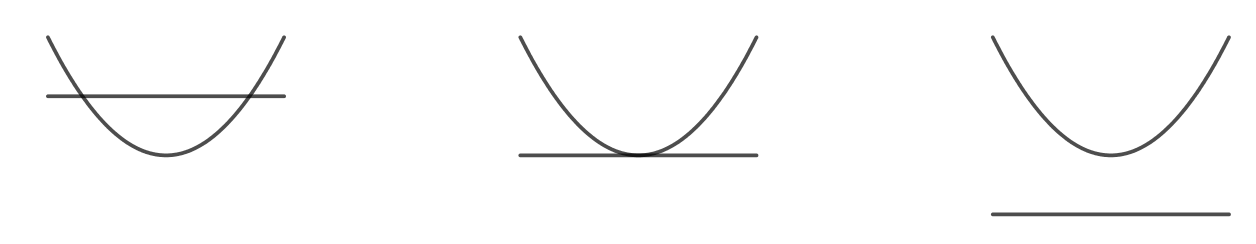}
\end{center}
\vspace{-0.3in}
\caption{Three instances of a conic meeting a line.}
\label{f2}
\end{figure}

Let's see where the circle with center $(x_0,y_0)$ and radius $r$ meets the line at infinity. The homogenization of the defining polynomial of the circle is $(x-x_0z)^2 + (y-y_0z)^2 - r^2z^2$ and setting $z=0$ we see that the circle meets the line at infinity at points $[x:y:0]$ where $x^2+y^2 = 0$. There are two such points, $I=[1:i:0]$ and $J=[1:-i:0]$, and they are independent of the radius and center of the circle. Everyone knows that circles are special kinds of conics, but now we know that they can be characterized as the irreducible conics that meet the line at infinity at the special points $I$ and $J$!   

\noindent{\bf Projective transformations.} In 1872 Felix Klein announced his Erlangen program for geometry, classifying geometries by the type of transformations that act on the underlying space and properties invariant under those maps. Projective space $\mathbb{P}^2$ admits \emph{projective transformations}, multiplication of points by an invertible $3 \times 3$ matrix $M \in GL_3(\mathbb{C})$. Such maps are well-defined since $M (\lambda \mathbf{v}) = \lambda (M \mathbf{v})$ and they send any collection of collinear points to new collinear points since multiplication by an invertible matrix preserves linear dependence. In fact, any four points in general position  in $\mathbb{P}^2$ -- no three collinear -- can be sent to any other four points in general position by a projective transformation. To see this, we first note that if the vectors $\mathbf{v_1}$, $\mathbf{v_2}$, $\mathbf{v_3}$, and $\mathbf{v_4} = a\mathbf{v_1} + b\mathbf{v_2}  +c\mathbf{v_3}$ are representatives for the homogeneous coordinates of four points in general position, then multiplying by the matrix $M$ whose columns are $a\mathbf{v_1}$, $b\mathbf{v_2}$, and $c\mathbf{v_3}$ sends the three standard basis vectors $\mathbf{e_1}$, $\mathbf{e_2}$, $\mathbf{e_3}$ and $\mathbf{e_1} + \mathbf{e_2} + \mathbf{e_3}$ to the four given points. If multiplying by the matrix $N$ sends the four points $\mathbf{e_1}$, $\mathbf{e_2}$, $\mathbf{e_3}$ and  $\mathbf{e_1} + \mathbf{e_2} + \mathbf{e_3}$ to four new points, then multiplying by $NM^{-1}$ sends $\mathbf{v_1}$, $\mathbf{v_2}$, $\mathbf{v_3}$, and $\mathbf{v_4}$ to these four new points too.  

Projection from a point gives a geometric example of a projective transformation from one line to another. Given a point $O$ and two lines $\ell$ and $\ell'$ not passing through $O$ then the projection map $\pi: \ell \rightarrow \ell'$ given by $\pi(P) = OP \cap \ell'$, and illustrated in Figure \ref{fig:proj}, is obtained by crossing the normal vector to $\ell'$ with the vector obtained from the cross product of the homogeneous coordinates of $O$ and $P$. The reader can check that this map can also be represented as multiplication by a matrix depending on $O$ and $\ell'$. 

\begin{figure}[h!t]
\begin{center}
\includegraphics[scale=0.16]{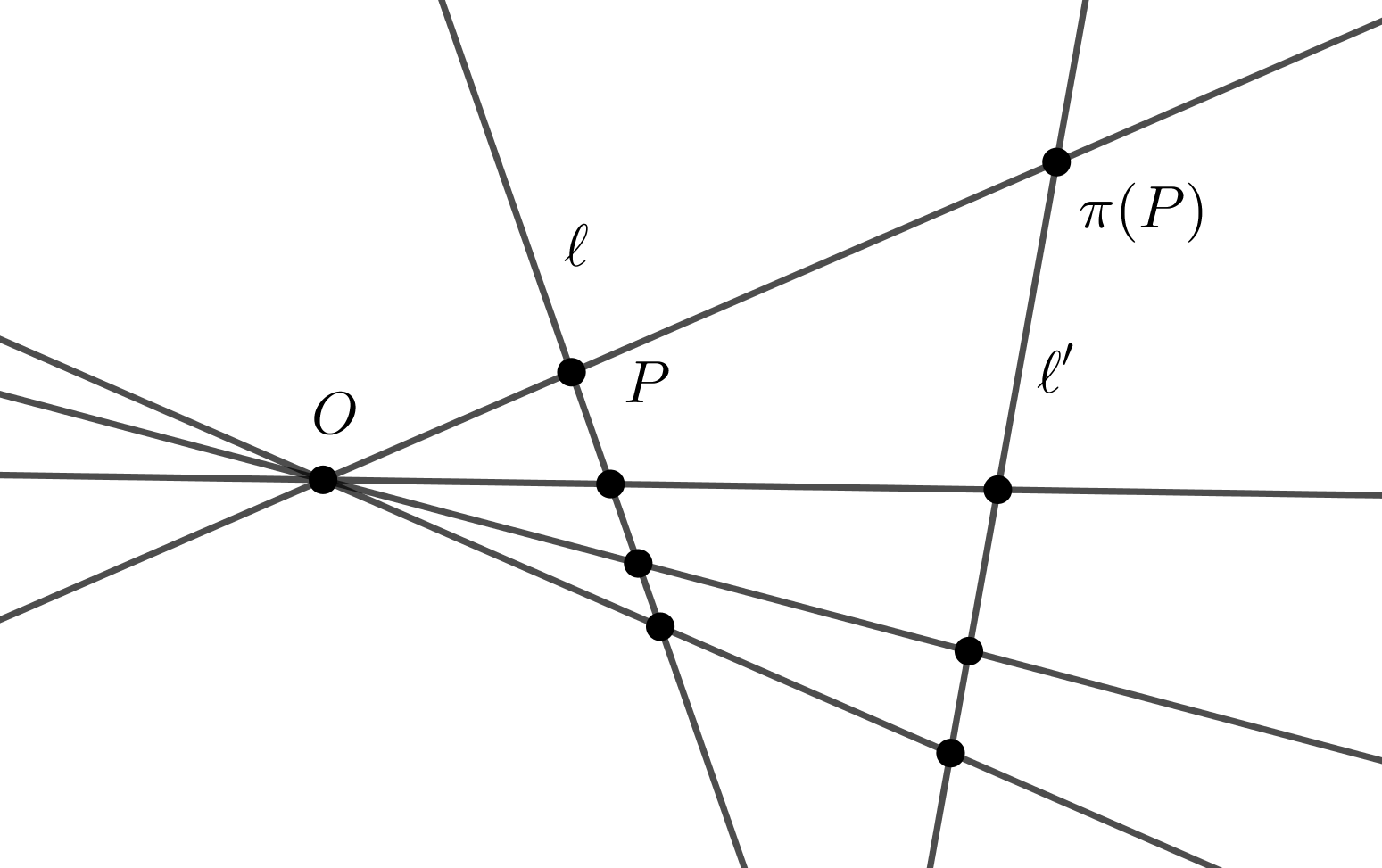}
\end{center}
\vspace{-0.3in}
\caption{The projection $\pi$ from $\ell$ to $\ell'$ through the point $O$.}
\label{fig:proj}
\end{figure}

Returning to the Erlangen program, multiplying by a matrix can be interpreted as changing the basis of the underlying space $\mathbb{C}^3$, so projective geometry is mainly concerned with properties that are independent of the choice of coordinates on our space. Incidence relations are geometric properties that do not depend on our coordinate system. The property of six points lying on a conic is also independent of the coordinates since Pascal's Theorem reduces the question of whether six points lie on a conic to an incidence statement. However, projective transformations do not preserve distances or angles, so the Erlangen program suggests that properties depending on distances and angles are not proper objects of study in projective geometry. Our point of view in this paper carefully blends projective and Euclidean geometry, using distances to prove results in Euclidean geometry and interpreting those results in terms of projective objects. For instance, to check whether six points lie on a conic we will apply a change of coordinates to move two of those points to $I$ and $J$, which reduces the problem to checking whether the remaining transformed points lie on a \emph{circle}, which we can check by measuring distances. This viewpoint seems in conflict with Klein's original Erlangen program but is consistent with his response to later developments \cite{Klein}.  Jacques-Salomon Hadamard \cite{Hadamard} wrote “It has been written\footnote{Apparently, Hadamard was paraphrasing  Paul Painlev\'{e} \cite{Painleve}, the French mathematician and statesman who served twice as Minister of War and twice as Prime Minister of France.} that the shortest and best way between two truths of the real domain often passes through the imaginary one.” Similarly, the best way between two truths in Euclidean geometry may pass through the projective domain.

\noindent{\bf Cross Ratios.} The cross ratio $(A,B; C,D)$ of four collinear points $A$, $B$, $C$ and $D$ is the quantity $$(A,B; C,D) = \frac{|AC||BD|}{|AD||BC|}.$$ 
Alexander Jones \cite{Jones} claims that Pappus of Alexandria already knew about cross ratios in 340 A.D but it was Carnot who introduced the use of oriented distances in the cross ratio. We fix a direction on the line to be positive and let $|AB|$ be the signed displacement from $A$ to $B$. 
If $O$ is any point off the line, the cross ratio $(A,B; C,D)$ equals $[OAC][OBD]/[OAD][OBC]$, where $[ABC]$ is the determinant of the $3 \times 3$ matrix whose columns are the homogeneous coordinates of the points $A$, $B$, and $C$.  Going forward, we will implicitly assume this interpretation of the distance ratios, allowing us to consistently orient all of our lines. With this interpretation, if $B$ is the point at infinity on the line then $[OBD]=[OBC]$ and $(A,B;C,D) = |AC|/|AD|$ is a quotient of directed distances measured from the common point $A$. The cross ratio is invariant under projective transformation since if we multiply by a matrix $M$ then the determinant of the matrix $[(MA)(MB)(MC)]$ of transformed points equals $\text{det}(M)[ABC]$ and then the determinant $\text{det}(M)$ cancels from all terms in the fraction defining $(MA,MB;MC,MD)$, yielding $(A,B;C,D)$.  

\end{section}

\begin{section}{Geometric Constructions}
\label{GCsection}

In this section we deal with some problems in constructive synthetic geometry: given enough information to uniquely define two geometric objects, we show how to construct their  intersection with a straightedge.

\noindent {\bf Construction 1: Find the second point of intersection of a line with a conic.} Given five points $P_1, \ldots, P_5$ lying on a unique conic $\mathcal{C}$ (that is, no four are collinear) and another point $Q$, we show how to use a straightedge to locate the second point of intersection $R$ of the line $\ell = \overlines{P_1Q}$ with the conic $\mathcal{C}$. This follows immediately from Pascal's Theorem if we arrange for one of the six lines in Pascal's Theorem to be the line $\overlines{P_1Q}$. Construct the three collinear points $Q_1 = \overlines{P_1P_2} \cap \overlines{P_4P_5}$, $Q_2 = \overlines{P_3P_4} \cap \ell$ and $Q_3 = \overlines{P_2P_3} \cap \overlines{Q_1Q_2}$. Then the point $R = \overlines{P_5Q_3} \cap \ell$ lies on $\mathcal{C}$ and $\ell$. This is illustrated in Figure \ref{int31} below. 

\begin{figure}[h!t]
\begin{center}
\includegraphics[scale=0.16]{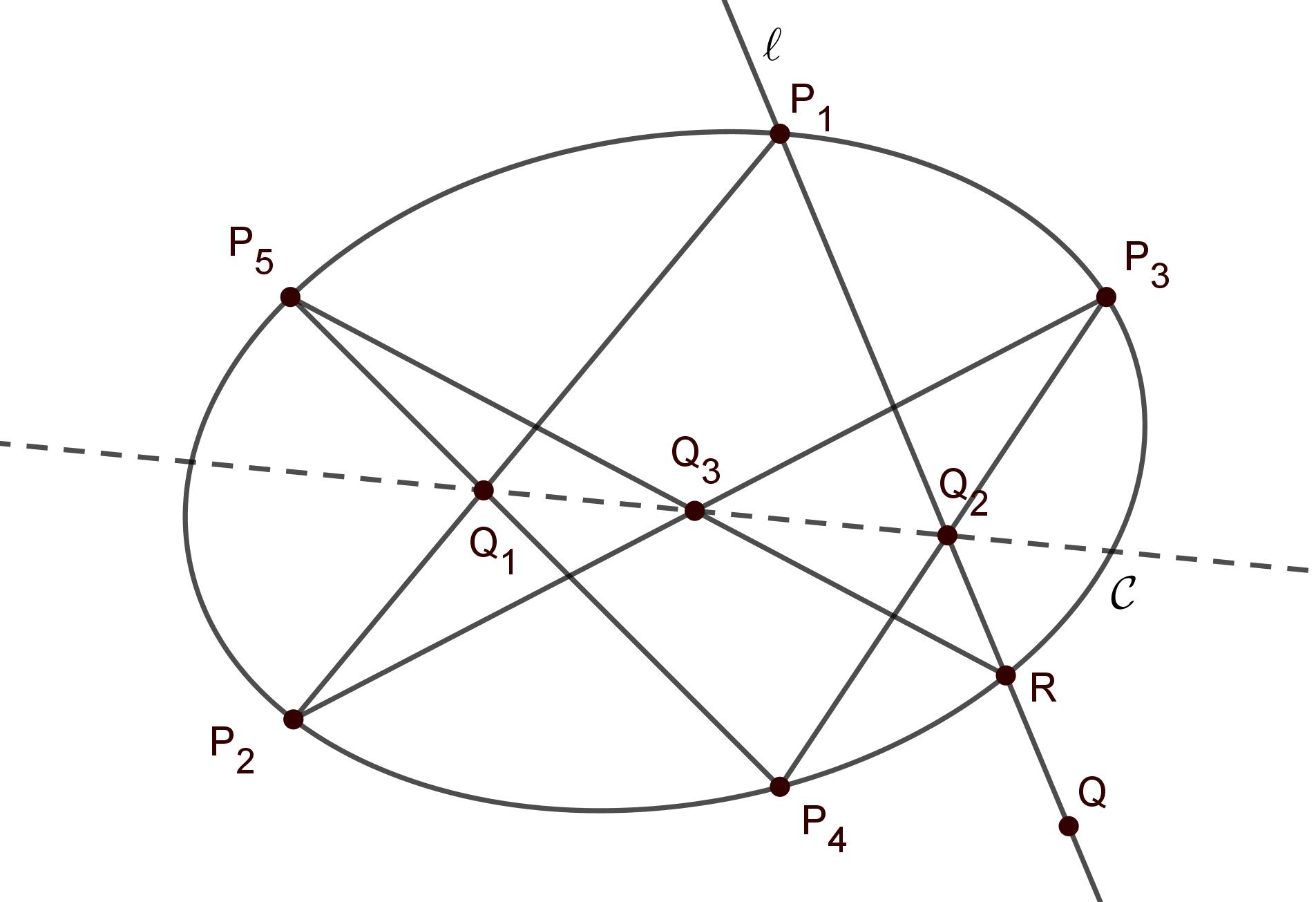}
\end{center}
\vspace{-0.3in}
\caption{Using Pascal's Theorem to find the fourth point of intersection of $\mathcal{C}$ with $\ell$.}
\label{int31}
\end{figure}

\noindent {\bf Construction 2: Find the fourth point of intersection of two conics.} Suppose that two conics $\mathcal{C}_1$ and $\mathcal{C}_2$ share the points $E_1$, $E_2$, and $E_3$ and that 
$\mathcal{C}_1$ also passes through points $P_1=[x_1:y_1:z_1]$ and $P_2=[x_2:y_2:z_2]$ and $\mathcal{C}_2$ passes through points 
$P_3=[x_3:y_3:z_3]$ and $P_4=[x_4:y_4:z_4]$. Further, suppose that all points are in general position so that 
no three are collinear. We show how to locate the fourth point of intersection $R$ of $\mathcal{C}_1$ with $\mathcal{C}_2$ using a straightedge. 

After a projective change of coordinates, we may assume that $E_1 = [1:0:0]$, $E_2 = 
[0:1:0]$ and $E_3 = [0:0:1]$. We choose to write the equation of the conic $\mathcal{C}_1$ as $a_1yz + a_2xz -a_3xy + \alpha_1 
x^2 + \alpha_2 y^2 + \alpha_3 z^2 = 0$ and the equation of $\mathcal{C}_2$ as $b_1yz + b_2xz -b_3xy + \beta_1 x^2 + 
\beta_2 y^2 + \beta_3 z^2 = 0$. Since both conics pass through $E_1$, $E_2$ and $E_3$, we find that all the $
\alpha$-coefficients and all the $\beta$-coefficients are zero. 

Now we define two maps. Let $E = \{[x:y:z] \in \mathbb{P}^2: \; xyz = 0\}$ and define $\phi: \mathbb{P}^2 \setminus E \rightarrow \mathbb{P}^2 \setminus E$ that sends the point $[x:y:z]$ to $[yz:xz:-xy]$. It will be convenient to think of the target $\mathbb{P}^2$ as the dual $\mathbb{P}^2$ so that the image $\phi([x_0:y_0:z_0])$ is the line $y_0z_0x + x_0z_0y - x_0z_0z = 0$. The second map $\Phi$ identifies the set of conics through $E_1$, $E_2$, and $E_3$ as a copy of $\mathbb{P}^2$. The map $\Phi$ sends the conic with equation $axz+bxz-cxy=0$ to the point $[a:b:c] \in \mathbb{P}^2$. The map $\Phi$ is well-defined: any other equation of the conic is a nonzero scalar multiple of the original equation and its associated point $[\lambda a: \lambda b: \lambda c]$ equals $[a:b:c]$. The two maps are closely related: 
$$ \Phi(\mathcal{C}) \in \phi(P) \Leftrightarrow P \in \mathcal{C},$$
that is, the point $\Phi(\mathcal{C}) \in \mathbb{P}^2$ lies on the line $\phi(P)$ precisely when the point $P \in \mathbb{P}^2$ lies on the conic $\mathcal{C}$. This follows immediately from the definitions of the two maps, but the reader is encouraged to pause and check this fact. 

Now the point $\Phi(\mathcal{C}_1)$ lies on both lines $\phi(P_1)$ and $\phi(P_2)$ so $\Phi(\mathcal{C}_1) = \phi(P_1) \cap \phi(P_2)$. Similarly, $\Phi(\mathcal{C}_2)$ is the intersection of the two lines $\phi(P_3)$ and $\phi(P_4)$. In fact, we can \emph{construct} the lines $L_i=\phi(P_i)$ with a straightedge using the process illustrated in Figure \ref{b2} (left). We join the points 
$[0:y_i:z_i] = E_1P_i \cap E_2E_3$ and $[x_i:0:z_i] = E_2P_i \cap E_1E_3$, producing the line $L_i$ with 
equation $y_iz_i x + x_iz_i y - x_iy_i z = 0$. You can check that this is the correct equation using the fact that that the cross product of two points gives the coefficients of the line through 
the points and the cross product of the vectors of coefficients for the lines gives the homogeneous coordinates of their intersection.

\begin{figure}[h!t]
\begin{center}
\includegraphics[scale=0.2]{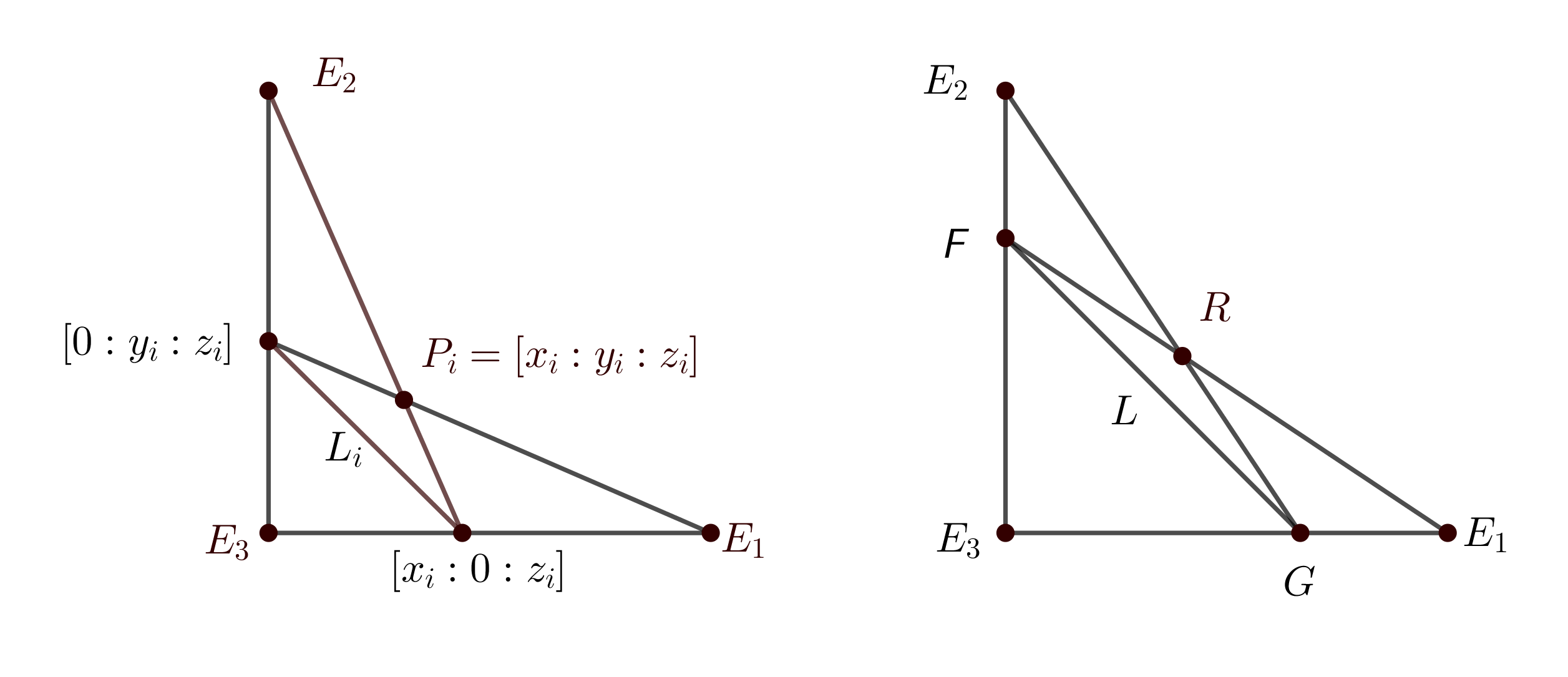}
\end{center}
\vspace{-0.4in}
\caption{Constructing $L_i$ from $P_i$ (left) and $R$ from $Q_1Q_2$ (right).}
\label{b2}
\end{figure}

So we can construct the two points $\Phi(\mathcal{C}_1)$ and $\Phi(\mathcal{C}_2)$ coming from the two conics. The line $L= \Phi(\mathcal{C}_1)\Phi(\mathcal{C}_2)$ joining these two points 
is the image $\phi(R)$ of a point $R$ that lies on both conics and is not $E_1$, $E_2$ or $E_3$. So $R$ is the fourth point of intersection we are looking for! We can reverse the process for constructing the image of $\phi$ to find the point $R$ with $\phi(R)=L$. The line $L$ 
meets the line $E_1E_3$ at $F$ and meets the line $E_2E_3$ at $G$. Then the lines $FE_1$ 
and $GE_2$ meet at $R$, constructing the fourth point of intersection of $\mathcal{C}_1$ and $\mathcal{C}_2$. The 
same construction works even when $E_1$, $E_2$ and $E_3$ do not lie in the special positions above since the 
projective transformation moving these points to the three given points of intersection of $\mathcal{C}_1$ and $\mathcal{C}_2$ 
preserves collinearities.

Let's pause our projective considerations to examine special results that hold for circles. In particular, we will make use of oriented distances. But don't worry, we will interpret these results in a projective way too. 

\noindent {\bf The Power of a Point Theorem.} Given a point $X$ and a circle, draw a line through $X$ that intersects the circle at $A$ and $B$. The {\em power of the point} $X$ with respect to this circle $\POP(X,\mathcal{C})$ is defined to be the product $|XA| \cdot |XB|$. The Power of a Point Theorem says that this quantity does not depend on the line drawn through $X$; if another line through $X$ meets the circle at $C$ and $D$ then $$ |XA| \cdot |XB| = |XC| \cdot |XD|. $$ 

\begin{figure}[h!t]
\begin{center}
\includegraphics[scale=0.2]{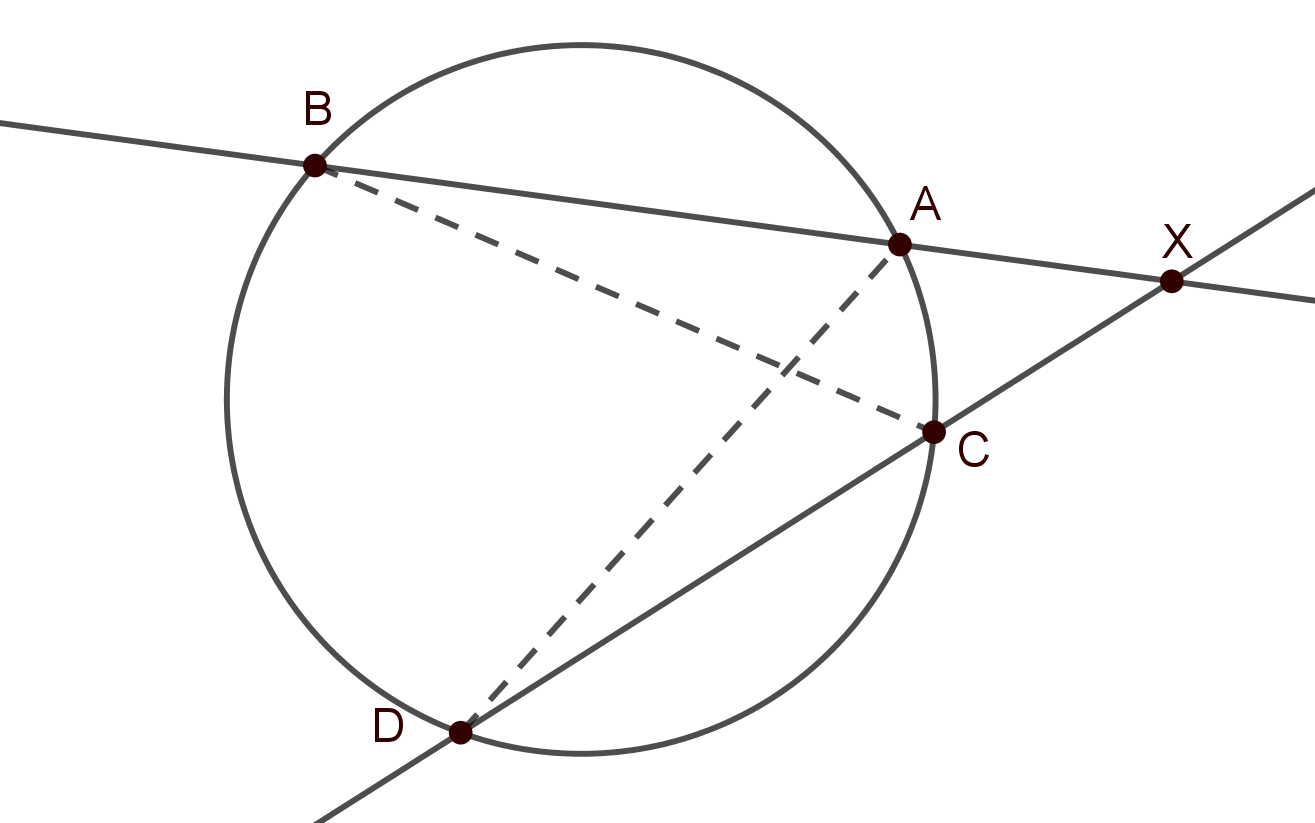}
\end{center}
\vspace{-0.3in}
\caption{The Power of a Point Theorem: $|XA| \cdot |XB| = |XC| \cdot |XD|$.}
\label{POPfig}
\end{figure}

The result follows easily from considerations about angles and similar triangles. Indeed, referencing Figure \ref{POPfig}, the two inscribed angles $\angle ADC$ and $\angle ABC$ are equal since both are half the central angle $\angle AOB$ supported at the center $O$ of the circle. Then the triangles $XAD$ and $XCB$ are similar so $$\frac{|XA|}{|XD|} = \frac{|XC|}{|XB|}$$ and hence $|XA| \cdot |XB| = |XC| \cdot |XD|.$

We will make use of the next result in our third construction.

\begin{lemma}\label{lemma:pc} Let $\mathcal{C}_1$ and $\mathcal{C}_2$ be irreducible conics meeting in four points $A$, $B$, $I'$ and $J'$. If a line through $A$ meets $\mathcal{C}_1$ and $\mathcal{C}_2$ at points $C$ and $D$, respectively, and a line through $B$ meets $\mathcal{C}_1$ and $\mathcal{C}_2$ at points $E$ and $F$, respectively, then the three lines $\overlines{CE}$, $\overlines{DF}$, and $\overlines{I'J'}$ are concurrent. In the special case where $I'=I$ and $J'=J$, the two conics are both circles and the two lines $CE$ and $DF$ meet on the line at infinity and hence are parallel.  \end{lemma}

\vspace{-0.15in}
\begin{figure}[h!t]
\begin{center}
\includegraphics[scale=0.2]{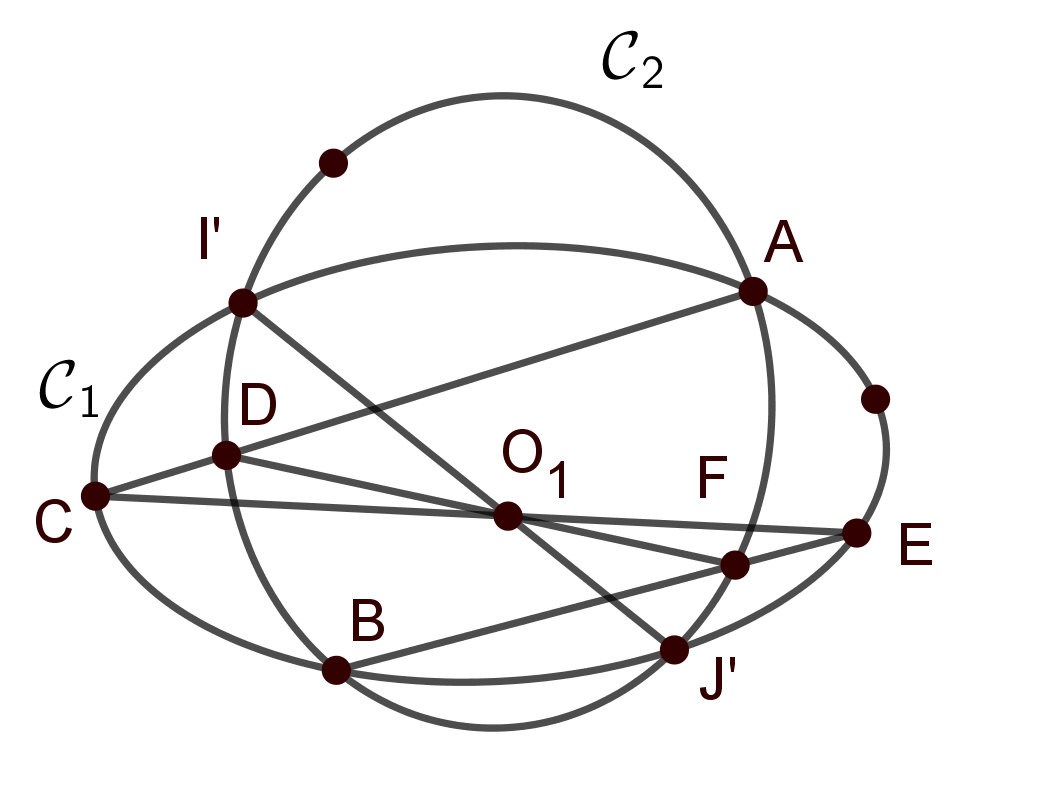}
\end{center}
\vspace{-0.3in}
\caption{The point $O_1$ lies on lines $CE$, $DF$ and $I'J'$.}
\label{parallelconics}
\end{figure}

\begin{proof}
Let $\mathcal{L}_A$ denote the line through $A$ and $\mathcal{L}_B$ denote the line through $B$.
Then Pascal's Theorem applied to the 6 points
$C, E, B, J', I', A$ of $\mathcal{C}_1$
tells us that the 3 points
$O_1 = \overlines{CE} \cap \overlines{I'J'}$,
$O_2 = \overlines{BE} \cap \overlines{AI'}=
\mathcal{L}_B \cap \overlines{AI'}$, and
$O_3 = \overlines{AC} \cap \overlines{BJ'}=
\mathcal{L}_A \cap \overlines{BJ'}$ are collinear. Similarly applying 
Pascal's Theorem to the 6 points
$D, F, B, J', I', A$ of $\mathcal{C}_2$
tells us that the 3 points
$O'_1 = \overlines{DF} \cap \overlines{I'J'}$,
$O'_2 = \overlines{BF} \cap \overlines{AI'}=
\mathcal{L}_B \cap \overlines{AI'}$, and
$O'_3 = \overlines{AD} \cap \overlines{BJ'}=
\mathcal{L}_A \cap \overlines{BJ'}$ are collinear.

Since $O_2=O'_2$ and $O_3=O'_3$ we see that
$O_1 = \overlines{O_2 O_3} \cap \overlines{I'J'} =
\overlines{O'_2 O'_3} \cap \overlines{I'J'} = O'_1$.
Hence $O_1 \in \overlines{CE} \cap \overlines{DF} \cap \overlines{I'J'}$.
\end{proof}

\noindent {\bf Construction 3: Given two conics sharing two known points, find the line through the remaining two points of intersection.} Consider the following situation. We are given 5 points $P_1$, $Q_1$, $R_1$ , $U$ and $V$, with no three collinear. These determine a unique irreducible conic $\mathcal{C}_1$ passing through the points. Further suppose we have 3 more points $P_2$, $Q_2$ and $R_2$ which together with $U$ and $V$ also determine a unique irreducible conic $\mathcal{C}_2$. Then $U$ and $V$ are two of the four points of intersection of $\mathcal{C}_1$ with $\mathcal{C}_2$. We now explain how to find the line through the other two points of the intersection of $\mathcal{C}_1$ with $\mathcal{C}_2$, which we call the radical axis of $\mathcal{C}_1$ and $\mathcal{C}_2$.

Pick a point $P$ and define
$C = \overlines{PU} \cap \mathcal{C}_1$
and
$D = \overlines{PU} \cap \mathcal{C}_2$.
Similarly we take 
$E = \overlines{PV} \cap \mathcal{C}_1$
and
$F = \overlines{PV} \cap \mathcal{C}_2$.
We can construct these 4 points $C,D,E,F$ using Construction 1. By Lemma~\ref{lemma:pc} the point
$O = \overlines{CE} \cap \overlines{DF}$ lies on the line through the other two points of intersection of $\mathcal{C}_1$ with $\mathcal{C}_2$.

Now we repeat the above with a new point $P'$ and
construct $C',D',E',F'$ and $O'$ from $P'$ in exactly the same manner.  Then the line $\overlines{OO'}$ is the desired
line through the other two points of intersection of $\mathcal{C}_1$ with $\mathcal{C}_2$.


\end{section}

\begin{section}{Ten Points on a Cubic}
\label{Section:10}

We are ready to determine if ten general points lie on a cubic curve. We first partition the ten points into two sets of five, $S_1$ and $S_2$. There is a conic $\mathcal{C}_i$ through the points in each $S_i$. Now swap two points of $S_1$ for two points of $S_2$ to make a second partition of the ten points into two sets of five, $T_1$ and $T_2$, so that $S_i$ and $T_i$ share three points in common. Again there is a conic $\mathcal{D}_i$ through the points in each $T_i$. The two unions of the pairs of conics form reducible degree-4 curves, $\mathcal{C} = \mathcal{C}_1 \cup \mathcal{C}_2$ and $\mathcal{D} = \mathcal{D}_1 \cup \mathcal{D}_2$. From B\'ezout's Theorem we know that $\mathcal{C}$ and $\mathcal{D}$ meet in 16 points, 10 of which were our original 10 points. The remaining six points will be referred to as the \emph{residual points}. The Cayley-Bacharach Theorem implies that the 10 original points lie on a cubic precisely when the six residual points lie on a conic. We will show how to locate the residual points and check that they lie on a conic using a straightedge.

To locate the residual points, note that the conics $\mathcal{C}_i$ and $\mathcal{D}_i$ share 3 known points, $S_i \cap T_i$, in common so we can construct their fourth point of intersection $P_i$ using Construction 2. As well, $\mathcal{C}_1$ and $\mathcal{D}_2$ share two known  points in common, $S_1 \cap T_2$, so we know that their two remaining points of intersection, $Q_1$ and $Q_2$, lie on a line $\mathcal{L}_Q$ that we can draw using Construction 3, even though we cannot construct the points $Q_1$ and $Q_2$ themselves. Similarly, $\mathcal{C}_2$ and $\mathcal{D}_1$ have two known points in common, $S_2 \cap T_1$, so the other two points of intersection, $R_1$ and $R_2$, lie on a line $\mathcal{L}_R$ that we can draw using Construction 3. The three lines $\mathcal{L}_P=P_1P_2$, $\mathcal{L}_Q$, and $\mathcal{L}_R$ form a triangle with vertices $P = \mathcal{L}_Q \cap \mathcal{L}_R$, $Q = \mathcal{L}_P \cap \mathcal{L}_R$, and $R=\mathcal{L}_P \cap \mathcal{L}_Q$. We would like to apply Carnot's Theorem to the six residual points $P_1$, $P_2$, $Q_1$, $Q_2$, $R_1$ and $R_2$, on this triangle to determine whether they lie on a conic. However, to use Carnot's Theorem we would need to compute terms like $|PQ_1|$ and $|PQ_2|$ and we do not know $Q_1$ and $Q_2$! Fortunately, we can use the Power of a Point Theorem to compute this quantity and achieve Galileo's dictum, ``make measurable what cannot be measured."
   
The two conics $\mathcal{C}_1$ and $\mathcal{D}_1$ intersect in four points, one of which is the residual point $P_1$. We choose one of the three other intersection points not on the line $PP_1$ and call it $G$. We label the other two intersection points $A$ and $B$. We draw the line $PG$, as illustrated in Figure \ref{thm-ill} (left). Construct $W$ as the second point of intersection of the line $\overlines{PG}$ with the conic $\mathcal{C}_1$ and construct $Z$ as the second point of intersection of the line $\overlines{PG}$ with the conic $\mathcal{D}_1$ using Construction 1. Noting that $P_1 \in \mathcal{L}_P$ lies on both $\mathcal{C}_1$ and $\mathcal{D}_1$, we can use Construction 1 to construct $X$ the second point of intersection of $\mathcal{L}_P$ with $\mathcal{C}_1$ and $Y$ the second point of intersection of $\mathcal{L}_P$ with $\mathcal{D}_1$. Now set $U = WX \cap \mathcal{L}_Q$ and $V = YZ \cap \mathcal{L}_R$.

\begin{figure}[h!t]
\includegraphics[scale=0.15]{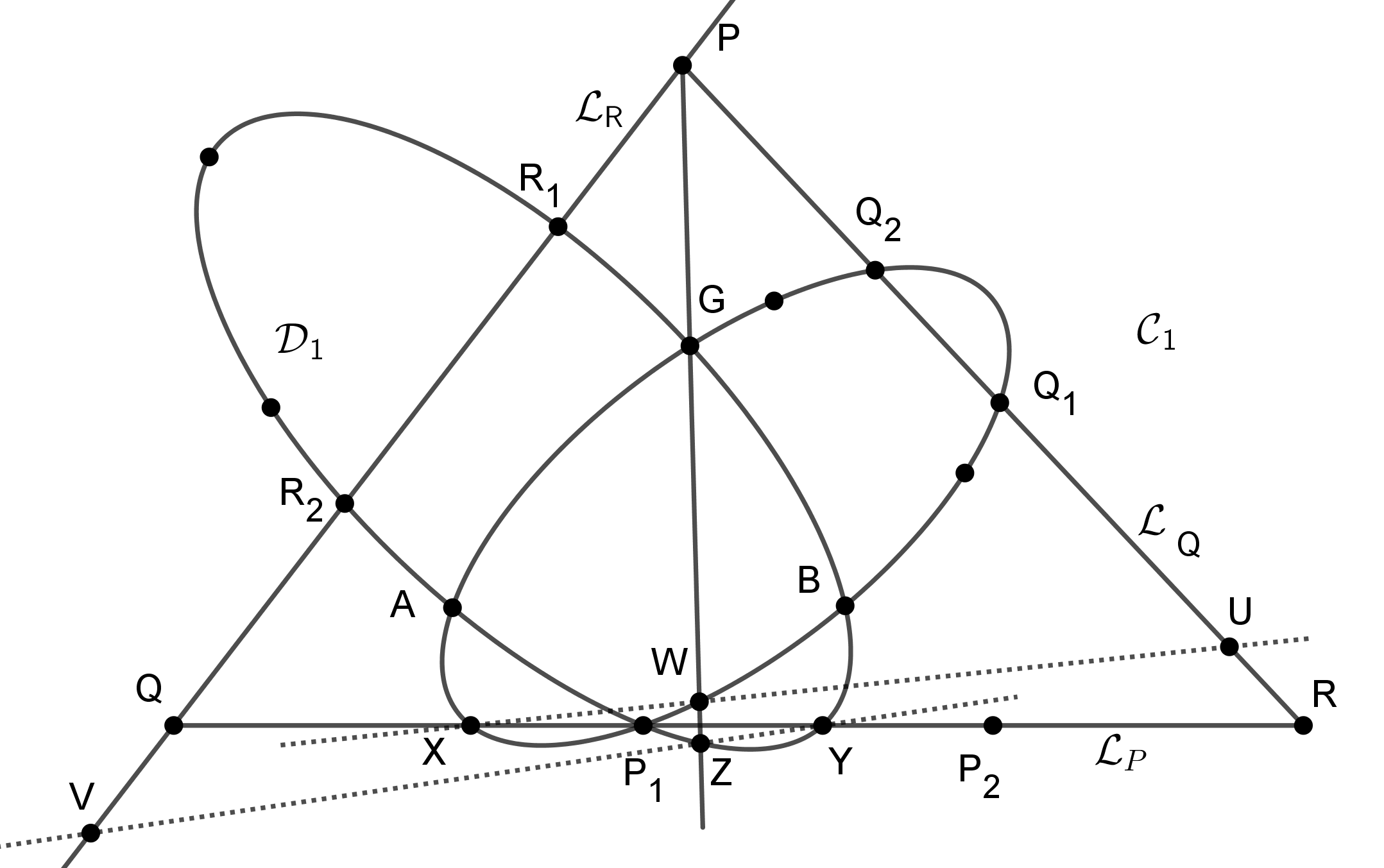}
\hspace{0.20in} \includegraphics[scale=0.075]{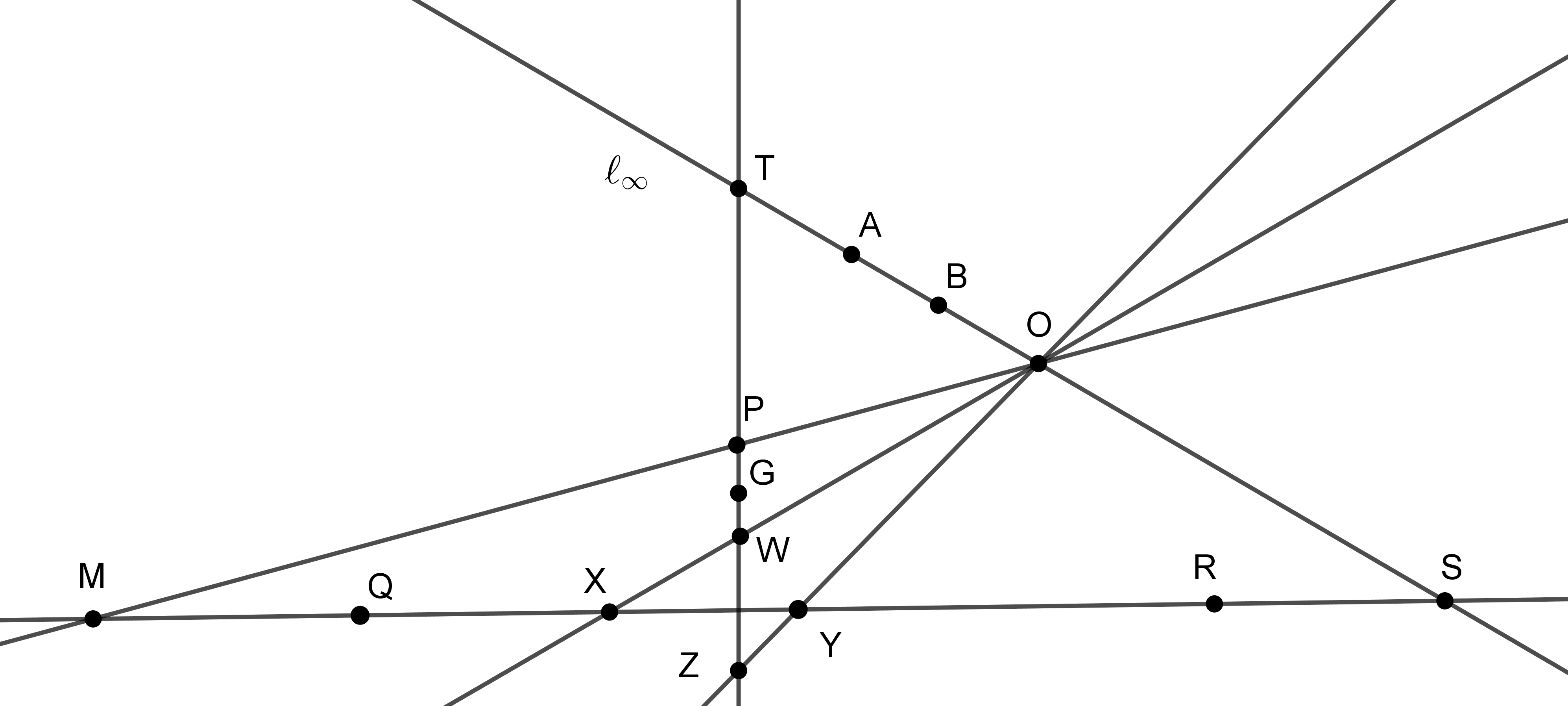} 
\caption{Two illustrations of the incidence relations.}
\label{thm-ill}
\end{figure}

\begin{theorem} \label{main-thm} The 10 general points lie on a cubic precisely when the six residual points lie on a conic and this happens if and only if the points $P_2$, $U$, and $V$ are collinear. \end{theorem}

Theorem \ref{main-thm} shows that we can check whether 10 points lie on a cubic curve using just an infinitely long complex straightedge on the projective plane. 

\begin{proof} As pictured in Figure \ref{thm-ill} (right), let $O = WX \cap YZ
$, $M = OP \cap RQ$, $S = AB \cap RQ$ and $T = AB \cap PG$. 

After applying a projective transformation we can assume  
$Q$ is $[0:0:1]$, $Y$ is $[1:0:1]$, $A$ is $I=[1:i:0]$ and $B$ is $J = [1:-i:0]$. The transformed conics $\mathcal{C}_1$ and $\mathcal{D}_1$ are circles since they pass 
through $I$ and $J$. The points $S$ and $T$ lie on the line at infinity $AB = IJ$. Moreover, the point $O$ now lies on the line at infinity since, by Lemma \ref{lemma:pc}, the two lines 
$WX$ and $YZ$ are parallel. Now consider the projection through $O$ mapping points on the line $
PG$ to points on the line $QR$ via $p \mapsto Op \cap 
QR$. This map sends $P$, $W$, and $Z$ to the points $
M$, $X$, and $Y$, respectively. Moreover, the point $T$ where the 
line $PG$ meets the line at infinity $AB=IJ$ is sent to the point 
$S$ where the line $QR$ meets the line at infinity since all three points $
S$, $T$ and $O$ are collinear.

 By Carnot's Theorem, the  six residual points $P_1$, $P_2$, $Q_1$, $Q_2$, $R_1$, and $R_2$, lie on a conic if and only if       

\begin{equation} \label{eq1}  \frac{|QP_1||QP_2| |RQ_1||RQ_2| |PR_1||PR_2|}{|RP_1||RP_2| |PQ_1||PQ_2| |QR_1||QR_2|}  = 1. \end{equation}

Now we reduce this condition using the Power of a Point several times. We have $$ \begin{array}{lll} 
|RQ_1||RQ_2| & = \POP(R,\mathcal{C}_1) & = |RP_1||RX| \\
|PR_1||PR_2| & = \POP(P,\mathcal{D}_1) & = |PG||PZ| \\
|PQ_1||PQ_2| & = \POP(P,\mathcal{C}_1) & = |PG||PW| \\
|QR_1||QR_2| & = \POP(Q,\mathcal{D}_1) & = |QP_1||QY|. \end{array}$$
Replacing terms in (\ref{eq1}), canceling three terms, and reordering the remaining products gives the equivalent condition
\begin{equation} \label{eq2} 
 \frac{|QP_2|}{|QY|} \frac{|RX|}{|RP_2|} \frac{|PZ|}{|PW|} = 1.
  \end{equation}
  
Each of the three ratios in the product in (\ref{eq2}) is the value of a cross ratio of four points. Recalling that $S$ is the point at infinity on $QR$ and $T$ is the point at infinity on $ZW$, we see that the six residual points lie on a conic precisely when 
\begin{equation} \label{eq3} 
 (Q,S;P_2,Y) (R,S; X,P_2) (P,T; Z,W)  =  1.
  \end{equation}

Using the above projection from $O$ to the line $\mathcal{L}_P = QR$ that sends $P$ to $M$, $T$ to $S$, $W$ to $X$, and $Z$ to $Y$, we can rewrite the term $(P,T; Z, W)$ as $(M,S; Y, X)$. Converting back to products of distance ratios and reordering we get an equivalent condition for the six residual points to lie on a conic
\begin{equation} \label{eq4} 
\frac{|P_2Q|}{|P_2R|} \frac{|XR|}{|XM|} \frac{|YM|}{|YQ|} = 1,
  \end{equation}
which can be expressed as a product of cross ratios $(P_2,S; Q,R) (X,S; R,M) (Y,S; M,Q)   =   1.$

Now apply the projection onto the line $PR$ via $O$ to the four points in the middle cross ratio. This maps $X$ to $U$, where $U = OX \cap PR$, and sends $S$ to the point at infinity $E=AB \cap PR = IJ \cap PR$ on the line $PR$. Also apply the projection from $O$ onto the line $PQ$ to the four points in the third cross ratio. This maps $Y$ to $V$, where $V = OY \cap PQ$, and sends $S$ to the point at infinity $F = AB \cap PQ = IJ \cap PQ$, on the line $PQ$. This produces the equivalent condition 
$$ \begin{array}{lllll} 
  & & 
(P_2,S; Q,R) (U,E; R,P) (V,F; P,Q)  & = &  \phantom{-}1 \\
& \Leftrightarrow & \frac{|P_2Q|}{|P_2R|} \frac{|UR|}{|UP|} \frac{|VP|}{|VQ|} & = & \phantom{-}1 \\
& \Leftrightarrow & \frac{|QP_2|}{|P_2R|} \frac{|RU|}{|UP|} \frac{|PV|}{|VQ|} & = & -1.
\end{array} 
$$
Using Menelaus's Theorem, this last condition can be interpreted as saying that the six points $P_1,$ $P_2,$ $Q_1,$ $Q_2,$ $R_1,$ and $R_2$ lie on a conic precisely when the three points $P_2$, $U$ and $V$ are the intersections of the extended edges of the triangle $PQR$ with a straight line. 

Applying the inverse transformation preserves the incidence relations: the original six residual points lie on a conic precisely when the three points $P_2$, $U=WX \cap PR$ and $V = YZ \cap PQ$ are collinear. This same condition checks whether the original 10 points lie on a cubic. \end{proof}  

We summarize our construction to check whether 10 points lie on a cubic. 

{\bf Construction 4: Checking whether 10 general points lie on a cubic.} \\
\vspace{-0.3in}
\begin{enumerate}
    \item Partition the 10 points into two sets $S_1$ and $S_2$, each with 5 points.
    \item Swap two points of $S_1$ with two points of $S_2$ to produce a second partition $T_1 \cup T_2$ of the 10 points with $| S_i \cap T_i | = 3.$ 
    \item Use Construction 2 to find the fourth points $P_i$ of intersection of the conic through $S_i$ with the conic through $T_i$ ($i = 1$ or $2$).
    \item Use Construction 3 to find the line $\mathcal{L}_Q$ through the two unknown points of intersection of the conic through $S_1$ and the conic through $T_2$. Use the same construction to find the line $\mathcal{L}_R$ through the two unknown points of intersection of the conic through $T_1$ and the conic through $S_2$. Use $\mathcal{L}_P = P_1P_2$, $\mathcal{L}_Q$ and $\mathcal{L}_R$ to construct $P = \mathcal{L}_Q \cap \mathcal{L}_R$, $Q = \mathcal{L}_P \cap \mathcal{L}_R$, and $R=\mathcal{L}_P \cap \mathcal{L}_Q$. 
    \item Take $G \in S_1 \cap T_1$ with $G \notin PP_1$. 
    Draw $PG$ and use Construction 1 to locate the intersection $W$ (resp. $Z$) of $PG$ with the conic through $S_1$ (resp. $T_1$).  
    \item Use Construction 1 to find the second point of intersection $X$ (resp. $Y$) of the line $\mathcal{L}_P = P_1P_2$ with the conic through $S_1$ (resp. $T_1$). 
    \item Take $U = XW \cap \mathcal{L}_Q$ and $V = YZ \cap \mathcal{L}_R$. 
    \item The 10 points lie on a cubic precisely when $P_2$, $U$ and $V$ are collinear.
\end{enumerate}

Our construction works for all sets of ten points off a set of measure zero. We can extend this construction to give an algorithm that will determine whether any set of 10 distinct points lies on a cubic curve. The extra work required is quite involved so we omit the details. The construction given here may fail when the 10 original points lie in special position. In particular, if our construction intersects two lines to create a point and the two lines are equal then the point is not well-defined. Similarly if we try to form the line through two points but the points are equal then the line is not well-defined.    
The degenerate configurations for which such problems occur possess additional geometric structure.   
We may exploit this extra structure using our straightedge to settle the question of whether the ten points lie on a cubic by ad-hoc means. As one example, if the 7 points of
$|S_1 \cup T_1|$ lie on a single conic then
$\mathcal{C}_1=\mathcal{D}_1$ and so Step~3 does not produce a well-defined point $P_1$. Besides just making a less fateful 
choice of $S_1$ and $T_1$, there is a simpler
way to proceed.  By
B\'ezout's Theorem, the cubic would have to be a union of the conic $\mathcal{C}_1$ and a line, so the ten points lie on a cubic precisely when the points off the conic $\mathcal{C}_1$ are collinear, which we can check using our straightedge. The other special cases are similar. 

\noindent {\bf Example.} We start with ten points:
$K_1=[ 0, 0, 1 ]$,  
$K_2=[ 6, -15, 1 ]$, 
$K_3=[ 1, 0, 1 ]$,  
$K_4=[ 2, 2, 1 ]$,  
$K_5=[ -\frac{5}{9}, \frac{8}{27}, 1 ]$,  
$K_6=[ 2, -3, 1 ]$,  
$K_7=[ \frac{1}{4}, -\frac{3}{8}, 1 ]$,  
$K_8=[ -1, 0, 1 ]$,  
$K_9=[ -1, -1, 1 ]$,  
$K_{10}=[ 6, 14, 1 ]$.  
Let $\mathcal{C}_1$ be the conic through $K_1,K_2,K_3,K_4,$ and $K_{5}$; $\mathcal{C}_2$ the conic through $K_6,K_7,K_8,K_9,$ and $K_{10}$; $\mathcal{D}_1$ the conic through $K_3,K_4,K_5,K_6,$ and $K_7$; and $\mathcal{D}_2$ the conic through $K_1,K_2,K_8,K_9,$ and $K_{10}$.  Conics $\mathcal{C}_1$ and $\mathcal{D}_1$ share three known points and meet in a fourth point 
$P_1 =[-10: 11: 1]$ and conics $\mathcal{C}_2$ and $\mathcal{D}_2$ also share three known points and meet in a fourth point $P_2 = [2: 5:1]$. Conics $\mathcal{C}_1$ and $\mathcal{D}_2$ share two known points and their radical axis is 
$ 10X -7Y + 18Z = 0$. Similarly, the radical axis of conics $\mathcal{C}_2$ and $\mathcal{D}_1$ is $2X - Y + 10Z = 0$. 
The two radical axes meet in $P = [ -13 : -16 : 1 ]$,  the first radical axis meets $P_1P_2$ in $R =[\frac{16}{9}: \frac{46}{9}: 1]$, and the second radical axis meets  $P_1P_2$ in $Q = [-\frac{8}{5} : \frac{34}{5} : 1]$. Conic $\mathcal{C}_1$ meets $\mathcal{L}_P$ in $P_1$ and $X = [\frac{12}{5} : \frac{24}{5} : 1]$ and conic $\mathcal{D}_1$ meets $\mathcal{L}_P$ in $P_1$ and $Y=[\frac{34}{11} : \frac{49}{11} : 1]$. The conics $\mathcal{C}_1$ and $\mathcal{D}_1$ meet at $G=K_3$, $A=K_4$,  $B=K_5$, and $P_1$. The line $PG$ meets $\mathcal{C}_1$ at $W =[ \frac{20}{13} : \frac{8}{13} : 1 ]$ and meets $\mathcal{D}_1$ at $Z =[ \frac{15}{22} : -\frac{4}{11} : 1 ]$. The line $WX$ meets the line $\mathcal{L}_Q$ at $U =[ \frac{11}{4} :  \frac{13}{2} : 1 ]$ and $YZ$ and $\mathcal{L}_R$ are parallel lines meeting at $V =[ \frac{1}{2} : 1 : 0 ]$ on the line at infinity. Now we can check that $P_2$, $U$ and $V$ all 
lie on the line $2X - Y + Z = 0$, which tells us that the 10 points we started with lie on a cubic curve. In fact, they all lie on the curve  $X^3-XZ^2-Y^2Z-YZ^2=0$.

The reader is warned that this example was meticulously constructed, involving very careful choices combined with a computer search among over a billion sets of points, to find an example involving rational coefficients expressible using small integers.  More typical examples generate coefficient explosion leading to coefficients for the constructed points requiring integers several tens or hundreds of digits long as numerators and denominators.

\end{section}

\begin{section}{A Binomial Proof} \label{sec:bp}

Having formulated Theorem \ref{main-thm}, we give a second, independent, proof of Theorem \ref{main-thm} found using an automated geometric theorem prover that we built.  
To explain how this works, we return to the bracket expressions $[ABC]$, representing the determinant of the $3\times 3$ matrix whose columns are the homogeneous coordinates of three points in $\mathbb{P}^2$. Consider the point $P$ that lies on the intersection of two lines $AB$ and $CD$. Then $P$ can be written as a linear combination of the points $A$ and $B$: $P = \alpha A + \beta B$. Since $P$, $C$ and $D$ are collinear, we have $[CDP]=0$ 
and substituting $\alpha A + \beta B$ for $P$ we find that $\alpha [CDA] + \beta [CDB] = 0$. It follows that $\alpha$ and $\beta$ can be chosen to be $[CDB]$ and $-[CDA]$, respectively; so $P = [CDB]A-[CDA]B$. Similarly, $P = [ABD]C - [ABC]D$. Equating the two expressions we have $[CDB]A-[CDA]B-[ABD]C = -[ABC]D$. Let $E$ be a fifth point and apply the operator taking the point $Q$ to $[QDE]$, giving (after some column interchanges) $$[BCD][ADE] - [ACD][BDE] - [ABD][CDE] = 0.$$ This expression is one of several quadratic relations among the brackets, called Grassmann-Pl\"{u}cker relations (see Richter-Gebert \cite{RG} for details). Now note that if $A$, $B$, and $D$ are collinear then $[ABD] = 0$ and so we get a \emph{binomial relation}, an equality of two bracket \emph{monomials}: $$ [BCD][ADE] = [ACD][BDE], $$ which we denote as $h(D,A,B,C,E)$. As well, this expression is equivalent to the collinearity of $A$, $B$ and $D$ if $[CDE] \neq 0$. 

The locus $\{X \in \mathbb{P}^2:  [ABX]=0\}$ is the line through $A$ and $B$. Similarly, $[ABX][CDX]=0$ represents a reducible conic, two lines whose union contains $A$, $B$, $C$ and $D$. Another reducible conic with the same property is $
[ACX][BDX]=0$. In fact, any conic through the four points has an equation of the form $\gamma [ABX][CDX] + 
\delta [ACX][BDX] =0$ for constants $\gamma$ and $\delta$. If we insist that the conic also pass through the 
point $E$ then $\gamma$ and $\delta$ can be taken to be $[ACE][BDE]$ and $-[ABE][CDE]$, respectively. That \
conic, $[ACE][BDE][ABX][CDX]-[ABE][CDE][ACX][BDX] = 0$, passes through a sixth point $F$ precisely when we 
have an equality of bracket monomials which we denote $c(A,B,C, D, E, F)$:  
$$ [ACE][ABF][CDF][BDE] = [ACF][ABE][CDE][BDF]. $$  

In the notation and setting of Theorem \ref{main-thm}, we have the following equalities from our collinearities: 

$$\begin{array}{lrll}  h( R_1,V, R_2, P_2 , Y) : &  [V, R_1, P_2][Y, R_2, R_1] & = & -[V,Y, R_1][ R_2, R_1, P_2]
\\      h( R_2, R_1, P, Q_1 , Z) : &  [ R_2, R_1, Q_1][Z, R_2, P] & = & -[Z, R_2, R_1][ R_2, Q_1, P]
\\      h( R_1, R_2, P, Q_2 , G) : &  [ R_2, R_1, Q_2][G, R_1, P] & = & -[G, R_2, R_1][ R_1, Q_2, P]
\\      h( P_1,Y, P_2, Q_2 ,  R_2) : &  [Y, Q_2, P_1][ R_2, P_2, P_1] & = & \phantom{-}[Y, R_2, P_1][ Q_2, P_2, P_1]
\\      h(Y, P_2, P_1,U , V) : &  [U,Y, P_2][V,Y, P_1] & = & \phantom{-}[V,Y, P_2][U,Y, P_1]
\\      h( P_1,X, P_2, Q_1 ,  R_1) : &  [X, Q_1, P_1][ R_1, P_2, P_1] & = & \phantom{-}[X, R_1, P_1][ Q_1, P_2, P_1]
\\      h(U, Q_2, P, P_2 , X) : &  [U, Q_2, P_2][U,X, P] & = & -[U,X, Q_2][U, P_2, P]
\\      h(Z,G, P, R_2 , Y) : &  [G,Z, R_2][Z,Y, P] & = & \phantom{-}[G,Z,Y][Z, R_2, P]
\\      h(G,W, P, P_1 ,  Q_1) : &  [G,W, P_1][G, Q_1, P] & = & \phantom{-}[G,W, Q_1][G, P_1, P]
\\      h(Y,V,Z, P ,  P_1) : &  [V,Y, P][Z,Y, P_1] & = & \phantom{-}[V,Y, P_1][Z,Y, P]
\\      h(V, R_1, P, P_2 , Y) : &  [V,Y, R_1][V, P_2, P] & = & -[V, R_1, P_2][V,Y, P]
\\      h( P_1,Y,X, Q_2 , U) : &  [U,Y, P_1][X, Q_2, P_1] & = & \phantom{-}[Y, Q_2, P_1][U,X, P_1]
\\      h( P_1,Y,X, R_1 , G) : &  [G,Y, P_1][X, R_1, P_1] & = & \phantom{-}[Y, R_1, P_1][G,X, P_1]
\\      h( Q_2,U, Q_1, P_2 , X) : &  [U,X, Q_2][ Q_2, Q_1, P_2] & = & -[U, Q_2, P_2][X, Q_2, Q_1]
\\      h( Q_2, Q_1, P, R_1 , W) : &  [W, Q_2, Q_1][ R_1, Q_2, P] & = & \phantom{-}[ R_1, Q_2, Q_1][W, Q_2, P]
\\      h( Q_1, Q_2, P, R_2 , G) : &  [G, Q_2, Q_1][ R_2, Q_1, P] & = & \phantom{-}[ R_2, Q_2, Q_1][G, Q_1, P]
\\      h(G,Z, P, P_1 ,  R_1) : &  [G,Z, R_1][G, P_1, P] & = & \phantom{-}[G,Z, P_1][G, R_1, P]
\\      h(W,G, P, Q_2 , X) : &  [G,X,W][W, Q_2, P] & = & \phantom{-}[G,W, Q_2][X,W, P]
\\      h(X,U,W, P ,  P_1) : &  [U,X, P_1][X,W, P] & = & \phantom{-}[U,X, P][X,W, P_1]

\end{array}$$

We also have equalities coming from six points lying on a conic: 

\begin{tabular}{rll}     {\makebox[3.1in][l]{$c(G,Z,Y, R_2, R_1, P_1)$:}}  &  & \\ \makebox[3.1in][r]{$[G,Y, R_1][G,Z, P_1][Y, R_2, P_1][Z, R_2, R_1]$} & = & $[G,Y, P_1][G,Z, R_1][Y, R_2, R_1][Z, R_2, P_1]$ \\
\makebox[3.1in][l]{$c(G, W,Q_2, P_1,X, Q_1):$} & &\\  \makebox[3.1in][r]{$[G,X, Q_2][G,W, Q_1][ Q_2, Q_1, P_1][X,W, P_1]$}  & = &  $[G, Q_2, Q_1][G,X,W][X, Q_2, P_1][W, Q_1, P_1]$  \\
\makebox[3.1in][l]{$c( R_1, Q_1,P_2, P_1,  R_2, Q_2):$} & & \\  \makebox[3.1in][r]{$[ R_2, R_1, P_2][ R_1, Q_2, Q_1][ Q_2, P_2, P_1][ R_2, Q_1, P_1]$}  & = &  { $[ R_1, Q_2, P_2][ R_2, R_1, Q_1][ R_2, P_2, P_1][ Q_2, Q_1, P_1]$ } 
\\ \makebox[3.1in][l]{$c(G, Z, R_1, P_1,Y, R_2):$} & & \\  \makebox[3.1in][r]{$[G, R_2, R_1][G,Z,Y][Y, R_1, P_1][Z, R_2, P_1]$} & = &  {$[G,Y, R_1][G,Z, R_2][ R_2, R_1, P_1][Z,Y, P_1]$ } 
\\ \makebox[3.1in][l]{$c(G, P_1, Q_2, Q_1, X,W):$} & & \\  \makebox[3.1in][r]{$[G,W, Q_2][G,X, P_1][X, Q_2, Q_1][W, Q_1, P_1]$} & = &  {$[G,X, Q_2][G,W, P_1][W, Q_2, Q_1][X, Q_1, P_1]$ }
\\ \makebox[3.1in][l]{$c( R_1, P_1, Q_2, Q_1, R_2, P_2):$} & & \\  \makebox[3.1in][r]{$[ R_1, Q_2, P_2][ R_2, R_1, P_1][ R_2, Q_2, Q_1][ Q_1, P_2, P_1]$}  & = &  {$[ R_2, R_1, Q_2][ R_1, P_2, P_1][ Q_2, Q_1, P_2][ R_2, Q_1, P_1]$ } \end{tabular}

Multiplying all the left-hand sides of these equalities together and multiplying all the right-hand sides together and canceling like terms leaves just a pair of brackets on each side: 
$$ [V, P_2, P][U,Y, P_2]  =  [V,Y, P_2][U, P_2, P] $$

This is precisely the equality $h( P_2,V,U, P , Y)$, so we find that $P_2$, $U$ and $V$ are collinear as long as $P,$ $P_2$ and $Y$ are not collinear (and all the brackets that we canceled are not zero, which requires a large collection of triples of points to be noncollinear). This is precisely the conclusion of Theorem \ref{main-thm}! This binomial proof was found using MATLAB \cite{MATLAB} to set up an integer programming problem, which we then passed to the optimization solver Gurobi \cite{gurobi} to solve. Our optimization problem involved 552 variables and 22,022 constraints, which Gurobi solved in just under 33 seconds on a five-year old laptop. We used the conclusion of Theorem \ref{main-thm} to set up the optimization problem, so this second proof only verifies the result in Theorem \ref{main-thm}; this method is not immediately applicable in the search for geometric results. 

We find this computational proof of Theorem \ref{main-thm} very amusing but it is hard to get any intuition for \emph{why} the result is true. This tendency of computer-assisted proofs to lead to results that cannot be easily explained to a human is one of the central problems with deep neural networks and other advanced tools in artificial intelligence today. 

\end{section}
\begin{section}{Extensions and Exercises}

We close with some fun problems, pointers to further reading, and comments about how this work connects to related topics. Two excellent books in geometry are Coxeter and Greitzer \cite{Coxeter} and Richter-Gebert \cite{RG}. The geometry chapter in Zeitz \cite{Zeitz} contains many challenging and enjoyable problems. The Cayley-Bacharach Theorem played a key role in our construction, but we never formally stated it. We refer the interested reader to the wonderful paper by Eisenbud, Green and Harris \cite{EGH}. The reader may also enjoy the enumerative study of conics in  Bashelor et al. \cite{BKT} and Traves \cite{Traves}, which contains another application of the Cayley-Bacharach Theorem to cubics. 


\noindent{\bf Problem 1.} Use the Power of a Point Theorem to prove Carnot's Theorem. A solution can be found at the wonderful Cut the Knot website. 

\noindent{\bf Problem 2.} Use Pascal's Theorem to construct the tangent line to a conic at a given point, given four additional points on the conic. As well, show that the tangent lines to a circle at $I$ and $J$ both pass through the center of the circle.  

\noindent{\bf Problem 3.} Define a transformation on homogeneous polynomials: if $M$ is a $3 \times 3$ matrix and $H(\mathbf{v}) = H(x,y,z)$ is a homogeneous polynomial of degree $d$, then define the transformed polynomial $M\cdot H$ to be $(M\cdot H)(\mathbf{v}) = H(M^{-1}\mathbf{v})$.  Check that $M\cdot H$ is a homogeneous degree-$d$ polynomial and that $H(\mathbf{v}) = 0$ precisely when $(M\cdot H)(M\mathbf{v}) = 0$, so changing the basis of $\mathbb{P}^2$ does not affect whether a collection of points lies on a degree-$d$ curve. 

\noindent{\bf Problem 4.} Four points $A,$ $B$, $C$, and $D$ whose cross ratio $(A,B; C,D)$ is $-1$ are said to be in \emph{harmonic position}. Let $x$ and $y$ be real numbers. Show that $$(-x,x;0,\infty) = (x,y; \frac{x+y}{2},\infty) = (-1,1;x,-x) = -1.$$ As well, show that permuting the four points in the cross ratio $(A,B; C,D)$ only produces six distinct values, rather than $4! = 24$ values. How are the six values related to one another? Which permutations fix the cross ratio?

\noindent{\bf Problem 5.} In 1857, Karl Georg Christian von Staudt introduced a way to view addition and multiplication in terms of cross ratios and incidence constructions. Fixing three points $0 = [0:0:1]$, $1 = [1:0:1]$ and $\infty = [1:0:0]$ and a fourth point $X = [x:0:1]$ on the $x$-axis  in $\mathbb{P}^2$ show that the cross ratio $(0, \infty; X, 1)$ equals $x$. If we follow the convention that $1/0 = \infty$ then $(0,\infty; \infty, 1) = \infty$ too, and we can identify the $x$-axis with the extended real numbers $\mathbb{R} \cup \{\infty\}$ via the cross ratio: the point $X$ corresponds to the value $(0, \infty; X, 1)$. The arithmetic operations of addition, subtraction, multiplication and division are defined on the extended real numbers and they induce corresponding operations on the projective line. Von Staudt's two key incidence structures are illustrated in Figure \ref{vs}, with some artistic license. For instance, given a second point $Y = [y:0:1]$ we construct a point $X+Y$ corresponding to the sum $x+y$ as follows. Pick a point $Q$ off the line containing $0$, $1$, $\infty$, $X$ and $Y$ and pick a second point $Q_1$ on the line joining $Q$ to $\infty$ (this line appears parallel to the $x$-axis in the finite part of the plane). Draw lines $\overlines{0Q}$ and $\overlines{XQ_1}$, meeting in point $P$. Draw lines $\overlines{P\infty}$ and $\overlines{YQ}$, meeting in point $R$. Then intersect the line $\overlines{RQ_1}$ with the original line $\overlines{0X}$ to obtain the point $X+Y$. Show that the resulting point $X+Y$ has the desired coordinates $[x+y:0:1]$, corresponding to the value $x+y$, irrespective of our choices for $Q$ and $Q_1$. As well, show that the incidence in Figure \ref{vs} (right) constructs the point $X\cdot Y=[xy:0:1]$. 

\begin{figure}[h!t]
\begin{center}
\includegraphics[scale=0.125]{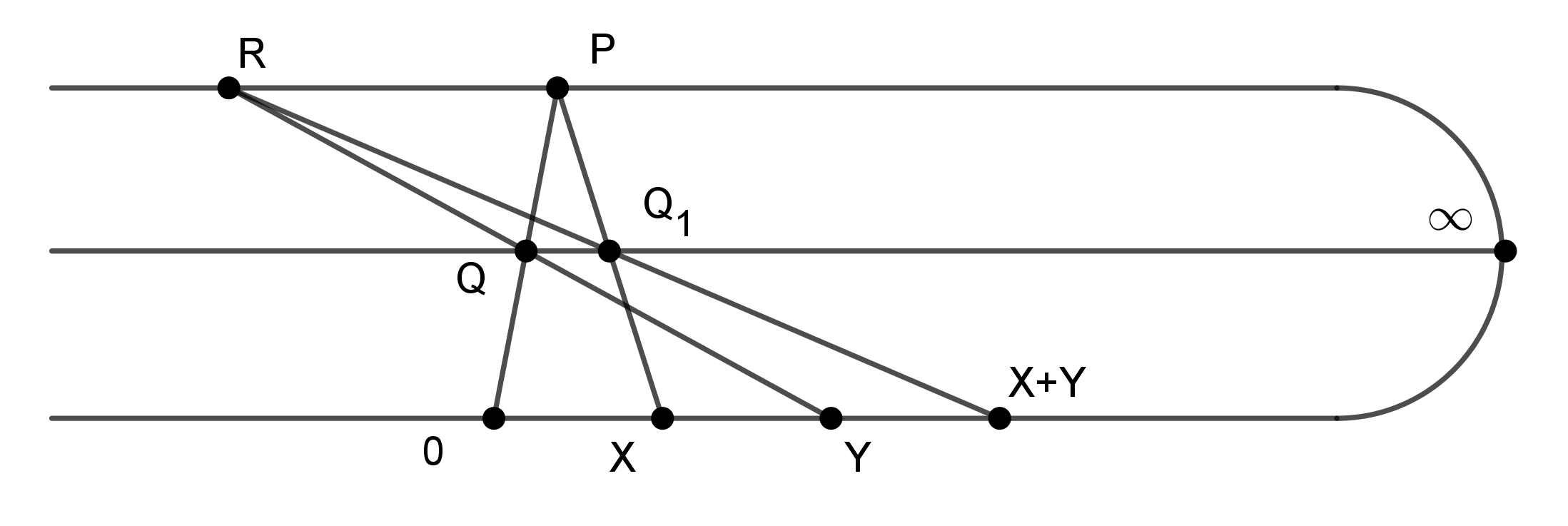}\hspace{0.25in} \includegraphics[scale=0.125]{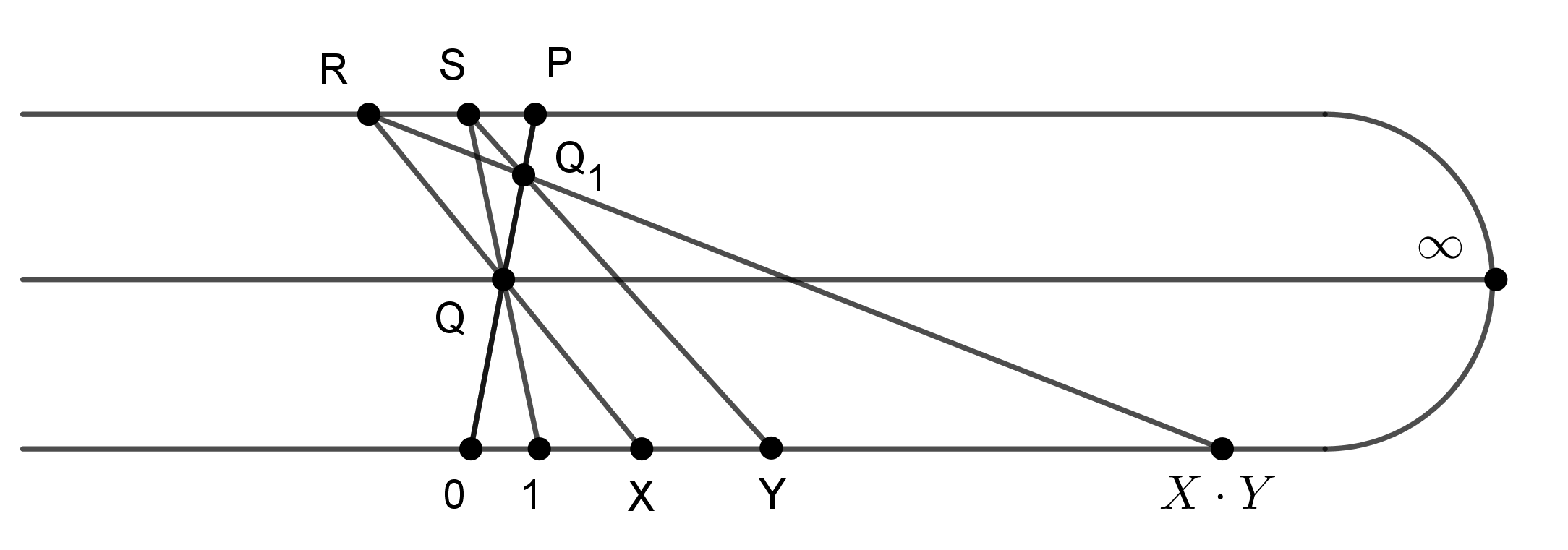}
\end{center}
\vspace{-0.3in}
\caption{Von Staudt's incidence structures that implement addition (left) and multiplication (right).}
\label{vs}
\end{figure}

\noindent{\bf Problem 6.} As in Construction 2, let $E_1 = [1:0:0],$ $E_2 = [0:1:0]$ and $E_3 = [0:0:1]$ and let $E = \{[x:y:z] \in \mathbb{P}^2: \; xyz = 0\}$. Show that the map $f: \mathbb{P}^2 \setminus E \rightarrow \mathbb{P}^2 \setminus E$ given by $f([x:y:z]) = [\frac{1}{x}: \frac{1}{y}: \frac{1}{z}]$ agrees with the map $\phi: \mathbb{P}^2 \setminus E \rightarrow \mathbb{P}^2 \setminus E$ given by $\phi([x:y:z]) = [yz: xz: -xy]$ and both are bijections. The map $\phi$ is actually defined on a slightly larger domain, $\mathbb{P}^2 \setminus \{ E_1, E_2, E_3\}$. What are the images of the lines $E_iE_j$ under the map $\phi$?  The map $\phi$ is an example of a Cremona transformation in algebraic geometry, a bijective map whose components are defined by polynomials from one dense open set of  $\mathbb{P}^n$ to another. 

\noindent{\bf Problem 7.} The $8 \Rightarrow 9$ theorem says that every plane cubic curve passing through 8 points of intersection of two cubic curves must pass through their ninth point of intersection as well. It follows that if these nine points of intersection are nine of our 10 points, then there is always a cubic through our ten points -- the result does not depend on the location of the tenth point! In this situation our construction to check whether ten points lie on a cubic is \emph{degenerate}. For instance, we need to draw a line through two points but these are the same point, or we need to intersect two lines, but these lines are the same line. 
Set up an example of this phenomenon and determine which kind of degeneracy occurs. We recommend working on this problem with a computer algebra system. 

Here is a matrix algebra approach to checking whether six points lie on a conic. For each point $P_i=[x_i:y_i:z_i]$ ($i\in \{1,\ldots,6\}$) we plug into the conic equation $ax^2+bxy+cxz + dy^2+eyz + fz^2 = 0$, producing a linear condition on the coefficients of the conic, $a, \ldots, f$. Together these six linear equations form a matrix equation $A\mathbf{v}=\mathbf{0}$, where $\mathbf{v} = [a,b,c,d,e,f]^T$. There is a conic through the six points precisely when this matrix $A$ has a nonzero nullspace. This occurs precisely when $\text{det}(A)=0$. This determinant condition can be written as a polynomial in the brackets $[P_iP_jP_k]$, each a determinant of a $3\times 3$ matrix. The determinant of $A$ is a degree-4 polynomial in the brackets with 720 summands. But the Grassmann-Pl\"ucker relations among the brackets allow us to reduce the number of summands. One can check that the polynomial evaluates to a scalar multiple of the expression $c(P_1,P_2,P_3,P_4,P_5,P_6)$ from the last section. 
The same approach can be used to check whether ten points lie on a cubic. In this case the matrix $A$ is a $10 \times 10$ matrix and in 1870 Reiss \cite{Reiss} wrote down an expression for $\text{det}(A)$ as a degree-10 polynomial in the brackets with just 20 terms. Suzanne Apel \cite{Apel} developed an algorithm to express a multiple of this polynomial by a monomial in the brackets as an expression in the \emph{Grassmann-Cayley algebra}, which can be interpreted as a massive incidence structure. Unfortunately, it is hard to get any intuition about why that structure implies the presence of a cubic through the ten points.  The reader might want to search for a short bracket polynomial expression for the determinant of the $15 \times 15$ matrix that determines whether 15 points lie on a degree-4 curve. This is an open research problem! 

We have focused on incidence structures in this paper and it is fun to consider \emph{extremal problems} with incidence structures. For instance, imagine a finite set of points for which every line through a pair of the points also contains a third point. One such collection consists of the nine points of inflection of a cubic curve. In fact, there are 12 lines passing through pairs (necessarily, triples) of the nine points, forming the beautiful incidence pattern depicted in Figure \ref{Hesse}. However, at least one of the points here has to have complex coordinates. The Sylvester-Gallai Theorem says that the only such point configuration with all real coordinates has all the points lying on a common line! 

\begin{figure}[h!t]
\begin{center}
\includegraphics[scale=0.2]{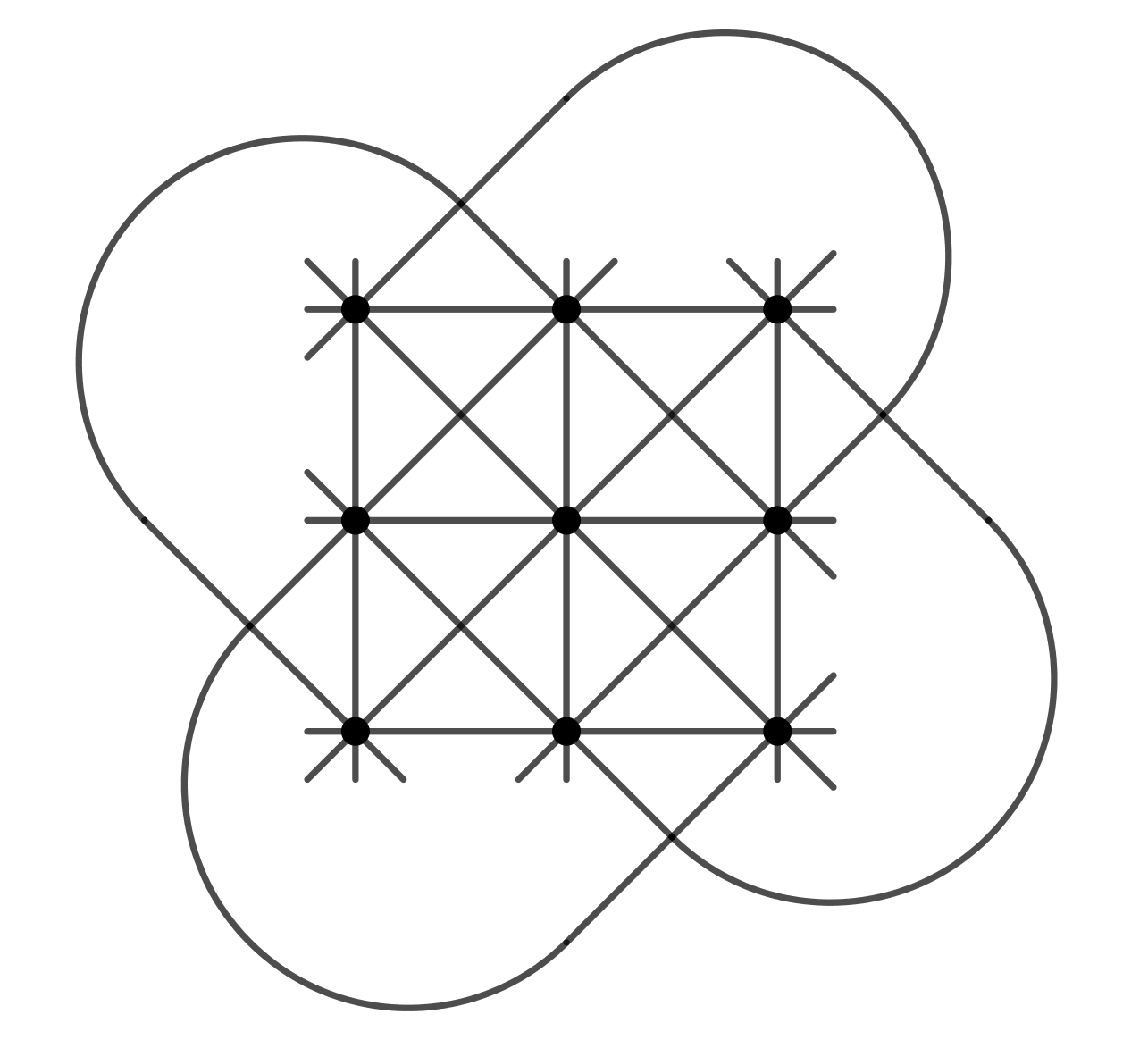}
\end{center}
\vspace{-0.3in}
\caption{The Hesse configuration: 12 complex lines through 9 points.}
\label{Hesse}
\end{figure}

\noindent{\bf Acknowledgments:} The authors thank Bernd Sturmfels for suggesting the problem to us and Mike Roth for helpful discussions. The computer algebra system MAGMA \cite{MAGMA} was also extremely helpful.   

\end{section}

\bibliography{10points}{}
\bibliographystyle{plain}

\end{document}